# WOODROOFE'S ONE-ARMED BANDIT PROBLEM REVISITED[1]

By Alexander Goldenshluger and Assaf Zeevi

*University of Haifa and Columbia University*

We consider the one-armed bandit problem of Woodroofe [*J. Amer. Statist. Assoc.* **74** (1979) 799–806], which involves sequential sampling from two populations: one whose characteristics are known, and one which depends on an unknown parameter and incorporates a covariate. The goal is to maximize cumulative expected reward. We study this problem in a minimax setting, and develop rate-optimal polices that involve suitable modifications of the myopic rule. It is shown that the regret, as well as the rate of sampling from the inferior population, can be finite or grow at various rates with the time horizon of the problem, depending on "local" properties of the covariate distribution. Proofs rely on martingale methods and information theoretic arguments.

## 1. Introduction.

*Background and motivation.* In his landmark paper, Robbins (1952) introduced an important class of sequential allocation problems, known collectively as multi-armed bandit problems. These models have since played central roles in areas such as statistics, operations research, engineering, computer science and economics. Berry and Fristedt (1985) and Gittins (1989) are standard references on the subject, and Lai (2001) and Cesa-Bianchi and Lugosi (2006) provide recent overviews of this voluminous literature.

The basic two-armed bandit problem can be described as follows. Consider two statistical populations. At each point in time, a single observation from one of the two populations can be taken, and a random "reward," governed by the properties of the sampled population, is realized. The objective is to devise a sampling policy that maximizes the expected cumulative (or

Received September 2008.
[1]Supported in part by BSF Grant 2006075.

discounted) reward over a designated time horizon. A particular case of interest arises when the probability distribution of one population is known a priori. This setting is often referred to as an *one-armed* bandit problem. Motivation for bandit problems can be found, among other examples, in clinical trials. Here the objective is to test two (or more) treatments and eventually allocate the one with greater efficacy to incoming patients; see Lai (2001) for further discussion and references therein.

The bulk of literature on multi-armed bandits assumes that the sampled populations are *homogeneous*. However, in many practical situations some additional information in the form of a covariate can be utilized for allocation purposes, and the reward distributions may depend on this covariate. For example, in the clinical trials context before deciding to assign a given patient to a treatment we can observe a covariate such as age or severity of disease. The one-armed bandit problem with covariates was first addressed in the pioneering work of Woodroofe (1979), who introduced and studied the following model.

*Woodroofe's one-armed bandit problem.* Let $(X_t, Y_t^{(0)}, Y_t^{(1)}, t \geq 1)$ denote a sequence of random vectors, where $X_t$ is a covariate value at stage $t$, and $Y_t^{(i)}$ is a potential reward from the arm $i = 0, 1$ that can be obtained at stage $t$. It is assumed that:

(i) the conditional distribution of $Y^{(0)}$ given $X$ is known, while the conditional distribution of $Y^{(1)}$ given $X$ depends on an unknown parameter $\theta$;

(ii) for any given value of the parameter $\theta$, $(X_t, Y_t^{(0)}, Y_t^{(1)})$ are independent and identically distributed (i.i.d.) copies of $(X, Y^{(0)}, Y^{(1)})$.

Suppose that $X_1, \dots, X_t, \dots$ are observed sequentially over time, and at each stage $t$ we can observe either $Y_t^{(0)}$ or $Y_t^{(1)}$, but not both. Let $\pi_t = 0$ or 1 according to whether $Y^{(0)}$ or $Y^{(1)}$ is observed at stage $t$; we will refer to this as sampling of arm 0 or arm 1, respectively. Then the objective is to develop a *sampling policy* $\pi = (\pi_t, t \geq 1)$ such that the expected value of the total geometrically discounted reward

$$V^* = \sum_{t=1}^{\infty} [\rho^{t-1}(\pi_t Y_t^{(1)} + (1 - \pi_t) Y_t^{(0)})]$$

is maximized. Here $\rho \in (0, 1)$ is a discount factor. By policy we mean a sequence $\pi = (\pi_t, t \geq 1)$ of random variables taking values in $\{0, 1\}$ such that $\pi_t$ depends on the observations collected up until time $t$.

Woodroofe (1979) considered the outlined problem in the Bayesian setting under the assumption that $Y^{(1)} = X - \theta + \varepsilon$, where $\varepsilon$ is a zero mean random variable with known distribution, independent of $X$. For a given prior distribution of $\theta$, Woodroofe (1979) provided a (nonconstructive) description of the optimal Bayesian policy. It was also shown that the *myopic policy*, which selects arm 1 when the value of $X$ is greater than the current estimate (posterior mean) of $\theta$, is asymptotically optimal as $\rho$ tends to 1. These results were later extended by Sarkar (1991) to a slightly more general setting.

*Summary of results.* In this paper we revisit Woodroofe's one-armed bandit problem with covariates, but with several notable distinctions.

We consider a minimax (non-Bayesian) framework with finite horizon $n$; see Section 2 for details. This is more in line with the formulation in Robbins (1952), and the seminal work of Lai and Robbins (1985) that developed asymptotically-optimal policies for traditional multi-armed bandits. The performance of a policy $\pi$ is measured relative to the *oracle policy* $\pi^* = (\pi^*_t, t \geq 1)$ that "knows" the unknown parameter $\theta$ and at each step, given the covariate value $X$, selects the arm with highest expected reward. In this context the *regret* and the *inferior sampling rate* are natural policy performance measures. The regret refers to the loss in the expected cumulative reward that stems from the use of a given policy relative to the oracle policy. The inferior sampling rate is the expected number of wrong arm selections prescribed by the policy. Assuming that the distribution of $(X, Y^{(1)})$ belongs to some natural class $\mathcal{P}$ of joint distributions, we measure performance of a policy $\pi$ by the maximal regret and inferior sampling rate over the class $\mathcal{P}$.

In this work we study minimax complexity of the one-armed bandit problem with covariates. We establish explicit nonasymptotic lower bounds on the minimax regret and inferior sampling rate (see Section 3.3, Theorems 3 and 4) and develop simple and intuitive *rate-optimal* policies which achieve these bounds in the sense of the order (see Section 3, Theorems 1 and 2).

Our work highlights a key property of the bandit problem with covariates: the performance of any policy depends critically on the behavior of the covariate distribution in the vicinity of the "decision boundary" $x = \theta$ (see Definition 1). This is akin to the Tsybakov margin condition that plays a pivotal role in nonparametric classification problems [see Mammen and Tsybakov (1999) and Tsybakov (2004a)]. In particular, depending on this margin condition, there are three distinct "regimes": one where it is possible to achieve a finite regret as $n \to \infty$; one where the regret grows like $\ln n$; and one where the regret grows like a fractional power of $n$ (see Remarks 3 and 4). These cases correspond to natural classes of distributions.

It is worth pointing out that the rate-optimal policies developed in this paper are not myopic. To that end, we were not able to prove that the myopic policy is rate optimal in our setting. This issue is discussed in Section 3.4, where our results are compared with those of Woodroofe (1979) and Sarkar (1991). The numerical results in Section 4 lend further credence to this point by illustrating the inferior performance of the myopic policy relative to the two policies proposed in this paper. See also further discussion in Section 4.

*Further related literature.* In contrast to the voluminous literature on traditional multi-armed bandits, the number of papers that address bandit problems with covariates is rather limited. We refer to Woodroofe (1982), Clayton (1989), Yang and Zhu (2002), Wang, Kulkarni and Poor (2005) and Goldenshluger and Zeevi (2008), where further references can be found.

The rest of the paper is organized as follows. Section 2 contains the problem formulation and the main definitions. In Section 3 we introduce two policies, establish upper bounds on their performance and derive lower bounds on the minimax regret and inferior sampling rate. Section 4 presents numerical results, and proofs of all results are given in Section 5.

## 2. Problem formulation.

*The model.* Assume that a sequence of i.i.d. random variables $X_1, X_2, \ldots$ with common distribution $P_X$ is observed sequentially over time. At each stage $t$, one can allocate the covariate $X_t$ to one of two response models obtaining random rewards $Y^{(0)}$ and $Y^{(1)}$, respectively. The random vectors $(X_t, Y_t^{(0)}, Y_t^{(1)})$ are i.i.d. copies of $(X, Y^{(0)}, Y^{(1)})$, and the conditional distribution of $Y^{(0)}$ given $X$ is known. Allocation of $X_t$ to the $i$th arm $(i = 0, 1)$ gives rise to a response (reward) $Y_t^{(i)}$ as follows:

$$Y_t^{(0)} = 0, \qquad Y_t^{(1)} = X_t - \theta + \varepsilon_t, \qquad t = 1, \ldots, n, \tag{2.1}$$

where $\theta$ is an unknown parameter, and $\{\varepsilon_t\}$ is a sequence of i.i.d. $\mathcal{N}(0, \sigma^2)$ random variables, independent of the sequence $(X_t, t \geq 1)$. As in Woodroofe (1979), the assumption $Y^{(0)} = 0$ does not restrict generality. Since the regret depends linearly on the observed rewards [see (2.2)], the reduction to $Y^{(0)} = 0$ is achieved by considering $Y = Y^{(1)} - \mathbb{E}(Y^{(0)}|X)$ instead of $Y^{(1)}$; thus, we always write $Y$ instead $Y^{(1)}$. Here and in what follows all random variables are assumed to be defined on the common probability space $(\Omega, \mathcal{F}, \mathbb{P})$, and $\mathbb{E}$ stands for the expectation with respect to $\mathbb{P}$.

*Policies and performance measures.* By a policy $\hat{\pi} = (\hat{\pi}_t, t \geq 1)$ we mean any sequence of random variables taking values in $\{0, 1\}$ such that $\hat{\pi}_t$ is $\mathcal{F}_{t-1}^+$-measurable; here $\mathcal{F}_{t-1}^+$ is the $\sigma$-field generated by the data collected up until time $t-1$, and by the current value of the covariate $X_t$, that is,

$$\mathcal{F}_t^+ := \sigma(X_1, \ldots, X_t, X_{t+1}; \hat{\pi}_1, \ldots, \hat{\pi}_t; \hat{\pi}_1 Y_1, \ldots, \hat{\pi}_t Y_t).$$

We also denote $\mathcal{F}_t := \sigma(X_1, \ldots, X_t; \hat{\pi}_1, \ldots, \hat{\pi}_t; \hat{\pi}_1 Y_1, \ldots, \hat{\pi}_t Y_t)$.

The quality of a policy $\hat{\pi} = (\hat{\pi}_t, t \geq 1)$ is measured relative to the performance of *the oracle* $\pi^* = (\pi^*_t, t \geq 1)$ which "knows" the value of the parameter $\theta$ a priori, and at each stage $t$ prescribes

$$\pi^*_t := \pi^*_t(\theta, X_t) = I\{X_t \geq \theta\}, \qquad t = 1, 2, \dots .$$

The *regret* $R_n(\hat{\pi}; \theta)$ is defined as the difference between the expected total rewards accumulated by the oracle $\pi^*$, and the expected total reward generated by $\hat{\pi}$ over a horizon $n$:

$$R_n(\hat{\pi}, \pi^*) := \mathbb{E}\sum_{t=1}^{n} \pi^*_t Y_t - \mathbb{E}\sum_{t=1}^{n} \hat{\pi}_t Y_t = \mathbb{E}\sum_{t=1}^{n} |X_t - \theta| I\{\hat{\pi}_t \neq \pi^*_t\}. \tag{2.2}$$

Another performance characteristic of a policy $\hat{\pi}$ is the *inferior sampling rate* defined by

$$S_n(\hat{\pi}, \pi^*) := \mathbb{E}[T_{\inf}(n)] = \sum_{t=1}^{n} \mathbb{P}\{\hat{\pi}_t \neq \pi^*_t\},$$

where $T_{\inf}(n) = \sum_{t=1}^{n} I\{\hat{\pi}_t \neq \pi^*_t\}$ is the total number of times the policy $\hat{\pi}$ sampled the inferior arm.

In this paper we adopt a minimax approach. Let $\mathcal{P}$ be a class of joint distributions $P_{X,Y}$ of $(X, Y)$; then the quality of a policy $\hat{\pi}$ will be measured by the maximal regret $R_n(\hat{\pi}, \mathcal{P}) = \sup_{P_{X,Y} \in \mathcal{P}} R_n(\hat{\pi}; \pi)$, and by the maximal inferior sampling rate $S_n(\hat{\pi}; \mathcal{P}) = \sup_{P_{X,Y} \in \mathcal{P}} S_n(\hat{\pi}, \pi^*)$. The minimax regret and the minimax inferior sampling rate are defined by

$$R^*_n(\mathcal{P}) := \inf_{\hat{\pi}} R_n(\hat{\pi}; \mathcal{P}), \qquad S^*_n(\mathcal{P}) := \inf_{\hat{\pi}} S_n(\hat{\pi}; \mathcal{P}),$$

where inf is taken over all possible policies. The policy $\hat{\pi}^*$ is said to be *rate optimal* with respect to the class $\mathcal{P}$ if

$$\limsup_{n\to\infty} \frac{R_n(\hat{\pi}^*; \mathcal{P})}{R^*_n(\mathcal{P})} < \infty \quad \text{and} \quad \limsup_{n\to\infty} \frac{S_n(\hat{\pi}^*; \mathcal{P})}{S^*_n(\mathcal{P})} < \infty.$$

In this paper we develop rate-optimal policies and study the behavior of the minimax inferior sampling rate and regret for some natural classes $\mathcal{P}$ of joint distributions $P_{X,Y}$.

*Classes of joint distributions $\mathcal{P}$.* It turns out that the complexity of the one-armed bandit problem with covariates (as measured by the minimax inferior sampling rate and the minimax regret) is essentially determined by the behavior of $P_X$ near the "decision boundary" $x = \theta$. This behavior can be quantified by a condition similar to the so-called *Tsybakov margin condition* in classification [cf. Mammen and Tsybakov (1999), Tsybakov (2004a)].

The joint distribution $P_{X,Y}$ can be described in terms of the conditional distribution $P_{Y|X}$ and the marginal distribution $P_X$. According to (2.1), the conditional distribution of $Y$ given $X$ is Gaussian with mean $X - \theta$ and variance $\sigma^2$, that is, $\mathcal{N}(x - \theta, \sigma^2)$.

DEFINITION 1. We say that $P_{X,Y} \in \mathcal{P}_\alpha(\theta)$ if $P_{Y|X=x} = \mathcal{N}(x-\theta, \sigma^2)$, and there exist constants $C_* > 0$, $\alpha > 0$, $x_0 \in (0, \frac{1}{2})$, $p \in (0,1)$ such that

$$P_X\{[\theta - x, \theta + x]\} \le C_* x^\alpha \qquad \forall x \in (0, x_0] \tag{2.3}$$

and

$$0 < p \le P_X\{[\theta, \infty)\} < 1. \tag{2.4}$$

Let $\Theta \subseteq \mathbb{R}^1$ and set

$$\mathcal{P}_\alpha(\Theta) = \bigcup_{\theta \in \Theta} \mathcal{P}_\alpha(\theta).$$

We write $\mathcal{P}_\alpha = \mathcal{P}_\alpha((-\infty, \infty))$.

Several remarks on the above definition are in order.

REMARK 1.

1. In what follows we omit in the notation explicit dependence of $\mathcal{P}_\alpha(\theta)$ on parameters $C_*, x_0$ and $p$. Without loss of generality, we suppose that the constants $C_*$, $x_0$ and $p$ are all the same for all classes $\mathcal{P}_\alpha(\theta)$. Note also that the parameters $C_*$, $x_0$, $\alpha$ and $p$ are related to each other; in what follows we assume that

$$p_1 := p - C_* x_0^\alpha > 0. \tag{2.5}$$

2. Condition (2.3) describes the behavior of $P_X$ near the "decision boundary" $x = \theta$. The most important and typical case is that of $\alpha = 1$: if $X$ has a density $f_X$ w.r.t. the Lebesgue measure separated away from zero in the vicinity of $x = \theta$, then $P_{X,Y} \in \mathcal{P}_1(\theta)$. If $f_X \sim |x-\theta|^{\alpha-1}$, $\alpha > 0$ for $x$ close to $\theta$, then $P_{X,Y} \in \mathcal{P}_\alpha(\theta)$. The case of $\alpha = \infty$ corresponds to a distribution of $X$ that assigns zero probability to the $x_0$-vicinity of $x = \theta$.
3. Condition (2.4) ensures that the oracle policy samples both from arm 0 and arm 1 with positive probability. In the absence of this condition, the problem is reduced to the setting with no covariates.

For the purpose of having a well-defined regret, we also consider the following restriction $\mathcal{P}'_\alpha(\Theta)$ of the class $\mathcal{P}_\alpha(\Theta)$.

DEFINITION 2. Assume that $\int |x| P_X(dx) < \infty$. For $\mu > 0$, let

$$\mathcal{P}'_\alpha(\theta) := \mathcal{P}_\alpha(\theta) \cap \left\{ P_{X,Y} : \int |x-\theta| P_X(dx) \le \mu \right\}$$

and $\mathcal{P}'_\alpha(\Theta) = \bigcup_{\theta \in \Theta} \mathcal{P}'_\alpha(\theta)$. We write $\mathcal{P}'_\alpha = \mathcal{P}'_\alpha((-\infty, \infty))$.

**3. Main results.** First we introduce some notation. For a policy $\pi = (\pi_t, t \geq 1)$, and $t = 1, 2, \ldots$, let $T_\pi(t) = \sum_{s=1}^t \pi_s$ denote the total number of times up until time $t$ that $\pi$ sampled from the arm 1. The estimator of $\theta$ based on the observations collected up until stage $t$ under the policy $\pi$ is defined by

$$\hat{\theta}_\pi(t) = \frac{1}{T_\pi(t)} \sum_{s=1}^{t} (X_s - Y_s)\pi_s. \tag{3.1}$$

3.1. *Nearly-myopic policy.* Consider the following policy $\hat{\pi}$. Set $\hat{\pi}_1 = 1$, and for $t = 1, 2, \ldots$, define

$$\hat{\pi}_{t+1} = I\left\{X_{t+1} \geq \hat{\theta}_{\hat{\pi}}(t) - \frac{\delta_t}{\sqrt{T_{\hat{\pi}}(t)}}\right\}, \tag{3.2}$$

where $\delta = (\delta_t, t \geq 1)$ is a sequence of positive real numbers to be specified. If $\delta_t \equiv 0$, for all $t$ in (3.2), then the corresponding policy is *myopic*, as it mimics the oracle policy $\pi^*$ by plugging in the current estimate of $\theta$.

The next theorem establishes nonasymptotic upper bounds on the maximal inferior sampling rate and the maximal regret of the policy $\hat{\pi}$.

THEOREM 1. *Let $\hat{\pi} = (\hat{\pi}_t, t \geq 1)$ denote the policy given by (3.1)–(3.2) and associated with $\delta_t = 2\sigma\sqrt{3 \ln t}$. Let $t_0 := \min\{t \in \{1, \ldots, n\} : x_0\sqrt{pt} \geq 8\sigma\sqrt{3 \ln t}\}$ and define*

$$t_\alpha := \min\{t \in \{1, \ldots, n\} : t \geq (8\sqrt{3}\sigma\sqrt{\ln t/p})^{4\alpha/(\alpha-2)}\} \qquad \forall \alpha > 2. \tag{3.3}$$

(i) *For all $\alpha > 0$,*

$$S_n(\hat{\pi}; \mathcal{P}_\alpha) \leq (t_0 \vee 2) + C_* \sum_{t=1}^{n} \left(8\sqrt{3}\sigma\sqrt{\frac{\ln t}{pt}}\right)^{\alpha} + K. \tag{3.4}$$

*Furthermore, if $\alpha > 2$, then*

$$S_n(\hat{\pi}; \mathcal{P}_\alpha) \leq (t_0 \vee t_\alpha \vee 2) + \frac{4C_*}{\alpha - 2} + K. \tag{3.5}$$

(ii) *For all $\alpha > 0$,*

$$R_n(\hat{\pi}; \mathcal{P}'_\alpha) \leq \mu[(t_0 \vee 2) + K] + C_* \sum_{t=1}^{n} \left(8\sqrt{3}\sigma\sqrt{\frac{\ln t}{pt}}\right)^{(\alpha+1)}. \tag{3.6}$$

*In addition, if $\alpha > 1$, then*

$$R_n(\hat{\pi}; \mathcal{P}'_\alpha) \leq \mu[(t_0 \vee t_{\alpha+1} \vee 2) + K] + \frac{4C_*}{\alpha - 1}. \tag{3.7}$$

The constant $K = K(p)$ appearing in (3.4)–(3.7) depends on $p$ only; its exact expression is given in the proof of the theorem.

REMARK 2. The bounds (3.4) and (3.6) are too conservative for large $\alpha$. In particular, they are not applicable to the case $\alpha = \infty$. On the other hand, (3.5) and (3.7) provide upper bounds on the inferior sampling rate and the regret in the important case of $\alpha = \infty$. In particular, $t_\infty \le c[(\sigma^2 p^{-1}) \ln(\sigma p^{-1/2})]^2$ for some constant $c$.

REMARK 3. The immediate consequence of Theorem 1 is that

$$S_n(\hat{\pi}; \mathcal{P}_\alpha) \le \begin{cases} C, & \alpha > 2, \\ C(\ln n)^2, & \alpha = 2, \\ Cn^{1-\alpha/2}(\ln n)^{\alpha/2}, & 0 < \alpha < 2, \end{cases} \tag{3.8}$$

$$R_n(\hat{\pi}; \mathcal{P}'_\alpha) \le \begin{cases} C, & \alpha > 1, \\ C(\ln n)^2, & \alpha = 1, \\ Cn^{(1-\alpha)/2}(\ln n)^{(1+\alpha)/2}, & 0 < \alpha < 1, \end{cases} \tag{3.9}$$

where $C$ depends on parameters of the class $\mathcal{P}_\alpha$ (resp., $\mathcal{P}'_\alpha$). Thus, the maximal inferior sampling rate of $\hat{\pi}$ is finite when $\alpha > 2$. Similarly, the maximal regret is finite when $\alpha > 1$. On the other hand, both the maximal inferior sampling rate and the maximal regret diverge to infinity when $\alpha \le 2$, and $\alpha \le 1$, respectively.

A natural question then is whether there exists a policy with slower growth rates for the inferior sampling rate (3.8) when $\alpha \le 2$, and for the regret (3.9) when $\alpha \le 1$.

3.2. *Forced sampling policy.* Let $q > 0$ be a design parameter to be specified. Define the sequence $\mathcal{T}_0 = (\tau_t : t \ge 1)$ of positive integers by $\tau_1 = 1$, $\tau_t = \lfloor \exp\{qt\} \rfloor$, $t \ge 2$. The number of elements $N_0(t)$ of the sequence $\mathcal{T}_0$ that are less than or equal to $t$ satisfies the following inequalities:

$$\frac{1}{q} \ln t - 1 \le N_0(t) \le \frac{1}{q} \ln(t+1).$$

Consider the subsequence $\mathcal{T}$ of $\mathcal{T}_0$ containing all nonequal elements of $\mathcal{T}_0$. It is easily seen that if

$$t \ge \nu := 1 + \frac{1}{q} \ln_+ \left( \frac{2}{e^q - 1} \right), \tag{3.10}$$

then $\tau_t - \tau_{t-1} \ge 1$; here $\ln_+(\cdot) = \max\{\ln(\cdot), 0\}$. Therefore, if $\tau_t \in \mathcal{T}_0$ and $t \ge \nu$, then also $\tau_t \in \mathcal{T}$. For all $t$, let $N(t) := \sum_{\tau \in \mathcal{T}} I\{\tau \le t\}$. Then $N(t) \le N_0(t)$

for all $t$, and

$$N(t) \le N_0(t) \le \frac{1}{q}\ln(t+1) \qquad \forall t, \tag{3.11}$$

$$N(t) \ge N_0(t) - N_0(\nu) \ge \frac{1}{q}\ln\left(\frac{t}{\nu+1}\right) - 1 \qquad \forall t > \nu. \tag{3.12}$$

Now we define the policy $\tilde{\pi} = (\tilde{\pi}_t, t \ge 1)$ in the following way: set

$$\tilde{\pi}_t = \begin{cases} 1, & t \in \mathcal{T}, \\ I\{X_t \ge \hat{\theta}_{\tilde{\pi}}(t-1)\}, & \text{otherwise}, \end{cases} \tag{3.13}$$

where $\hat{\theta}_\pi(t)$ is given in (3.1). Thus, policy $\tilde{\pi}$ incorporates forced sampling from arm 1 at time instants $\mathcal{T}$, and myopic action in between. Under the policy $\tilde{\pi}$, arm 1 is pulled at least $N(t)$ times up until time $t$.

Let $\nu$ be given in (3.10), and let

$$\nu_0 := \max\{\nu, \min(t : t \ge 2q^{-1}\ln(t+1))\}. \tag{3.14}$$

Theorem 2. *Let $\tilde{\pi} = (\tilde{\pi}_t, t \ge 1)$ be the policy defined in (3.13) and associated with parameter $q > 0$. Then for all $n \ge \nu_0$ and any class $\mathcal{P}_\alpha/\mathcal{P}'_\alpha$ satisfying $x_0^2 \ge 12q\sigma^2$, one has the following:*

(i) *For any $\alpha > 0$,*

$$S_n(\tilde{\pi}; \mathcal{P}_\alpha) \le C_*\varkappa_\alpha \sum_{t=1}^n \left(4\sigma\sqrt{\frac{2}{p_1 t}}\right)^\alpha + \frac{1}{q}\ln(n+1) + \frac{16\alpha\sigma^2}{x_0^2 p_1} + K_1, \tag{3.15}$$

*where $\varkappa_\alpha := [(\alpha/2)^{\alpha/2}(1-2^{-\alpha})^{-1} + \Gamma(\alpha/2)(2\ln 2)^{-1}]$, and $K_1$ is a constant depending on parameters of the class $\mathcal{P}_\alpha$, $\sigma^2$ and $q$, but independent of $\alpha$ and $C_*$. Furthermore, if $\alpha > 2$, then*

$$S_n(\tilde{\pi}; \mathcal{P}_\alpha) \le \frac{1}{q}\ln(n+1) + C_*\left\{\frac{1}{1-2^{2-\alpha}}\left(\frac{32\sigma^2}{x_0^2 p_1}\right)^{\alpha/(\alpha-2)} + K_2\right\} + K_1, \tag{3.16}$$

*where $K_2$ is an absolute constant.*

(ii) *For any $\alpha > 0$,*

$$\begin{aligned} R_n(\tilde{\pi}; \mathcal{P}'_\alpha) \le{}& C_*\varkappa_{\alpha+1} \sum_{t=1}^n \left(4\sigma\sqrt{\frac{2}{p_1 t}}\right)^{\alpha+1} \\ &+ \mu\left\{\frac{1}{q}\ln(n+1) + \frac{16(\alpha+1)\sigma^2}{x_0^2 p_1} + K_1\right\}. \end{aligned} \tag{3.17}$$

*Furthermore, if* $\alpha > 1$, *then*

$$\begin{aligned} R_n(\tilde{\pi}; \mathcal{P}'_\alpha) &\le \frac{\mu}{q}\ln(n+1) \\ &\quad + C_*\left\{\frac{1}{1-2^{1-\alpha}}\left(\frac{32\sigma^2}{x_0^2 p_1}\right)^{(\alpha+1)/(\alpha-1)} + K_2\right\} + \mu K_1. \end{aligned} \tag{3.18}$$

*The exact expressions for constants* $K_1$ *and* $K_2$ *are given in the proof of the theorem.*

REMARK 4.

1. Note that the statement of the theorem holds only for classes $\mathcal{P}_\alpha$ (resp., $\mathcal{P}'_\alpha$) such that $x_0^2 \ge 12q\sigma^2$. This is in contrast to the policy $\hat{\pi}$ for which the results of Theorem 1 hold for classes $\mathcal{P}_\alpha$ and $\mathcal{P}'_\alpha$ with arbitrary parameters. Thus, the smaller the design parameter $q$, the larger the class of joint distributions for which the theorem statement is valid. Note, however, that the regret and the inferior sampling rate grow as $q$ decreases.
2. An immediate consequence of the theorem is that

$$S_n(\tilde{\pi}; \mathcal{P}_\alpha) \le \begin{cases} C\ln n, & \alpha \ge 2, \\ Cn^{1-\alpha/2}, & 0 < \alpha < 2, \end{cases}$$

$$R_n(\tilde{\pi}; \mathcal{P}'_\alpha) \le \begin{cases} C\ln n, & \alpha \ge 1, \\ Cn^{(1-\alpha)/2}, & 0 < \alpha < 1. \end{cases}$$

3. Comparing the above bounds with (3.8) and (3.9), we conclude that the forced sampling policy $\tilde{\pi}$ is better than the nearly-myopic policy $\hat{\pi}$ in the zone of "small" $\alpha$ ($\alpha \le 2$ for the inferior sampling rate, $\alpha \le 1$ for the regret). However, the inferior sampling rate and the regret of $\tilde{\pi}$ grow at least logarithmically for all $\alpha$, so $\hat{\pi}$ is better in the zone of "large" $\alpha$. We were not able to develop a single policy that simultaneously shares properties of $\hat{\pi}$ for large $\alpha$ and $\tilde{\pi}$ for small $\alpha$.

3.3. *Lower bounds.* Theorem 1 shows that $S_n(\hat{\pi}; \mathcal{P}_\alpha)$ is finite when $\alpha > 2$; likewise, $R_n(\hat{\pi}; \mathcal{P}'_\alpha)$ is finite when $\alpha > 1$. The next theorem establishes lower bounds on $S_n(\hat{\pi}; \mathcal{P}_\alpha)$ and $R_n(\hat{\pi}; \mathcal{P}'_\alpha)$ when $\alpha \le 2$ and $\alpha \le 1$, respectively.

THEOREM 3. *For any policy* $\pi$ *and large enough* $n$, *one has*

$$S_n(\pi; \mathcal{P}_\alpha) \ge \frac{1}{8}\left(\frac{\alpha}{2e}\right)^{\alpha/2} C_*\sigma^\alpha n^{1-\alpha/2} \qquad \forall \alpha \in (0, 2],$$

$$R_n(\pi; \mathcal{P}'_\alpha) \ge \left(\frac{1}{8}\right)^{1+1/\alpha}\left(\frac{\alpha}{2e}\right)^{(\alpha+1)/2} \frac{C_*^{1+1/\alpha}\sigma^{\alpha+1} n^{(1-\alpha)/2}}{2\max\{(1/x_0), (2C_*)^{1/\alpha}\}} \qquad \forall \alpha \in (0, 1].$$

Theorem 3 shows that the forced sampling policy $\tilde{\pi}$ is rate optimal with respect to the inferior sampling rate when $\alpha \in (0,2)$ and with respect to the regret when $\alpha \in (0,1)$. At the same time, in the zone of small $\alpha$ the performance of the nearly-myopic policy $\hat{\pi}$ is worse than the minimax rate by a logarithmic factor (see Theorem 1). Note, however, that in the "boundary" cases $\alpha = 2$ and $\alpha = 1$, there is a logarithmic gap between the lower bounds of Theorem 3 and the upper bounds of Theorem 2. This raises the question whether the forced sampling policy is rate optimal with respect to the regret (inferior sampling rate) whenever $\alpha = 1$ ($\alpha = 2$). Next we show that, for a wide class of admissible policies, the performance of the forced sampling policy cannot be improved upon.

Let $\Pi$ denote the class of policies $\pi = (\pi_t, t \geq 1)$ of the form $\pi_t = I\{X_t \geq \hat{\gamma}_t\}$, where $\hat{\gamma}_t$ is an $\mathcal{F}_{t-1}$-measurable random variable. We note that the class $\Pi$ is sufficiently rich and include policies with forced sampling (set, e.g., $\hat{\gamma}_t = \pm\infty$). We have the following result.

THEOREM 4. *Let $\Theta \subset \mathbb{R}^1$ be a closed bounded interval; then for all $n$ large enough, one has*

$$\inf_{\pi\in\Pi} S_n(\pi; \mathcal{P}_2(\Theta)) \geq K_1\sigma^2 \ln n,$$

$$\inf_{\pi\in\Pi} R_n(\pi; \mathcal{P}_1'(\Theta)) \geq K_2\sigma^2 \ln n,$$

*where $K_1$, $K_2$ are absolute constants.*

Thus, Theorem 4 establishes that in the "boundary" case ($\alpha = 1$ and $\alpha = 2$, for the regret and inferior sampling rate, resp.) the forced sampling policy $\tilde{\pi}$ cannot be improved in the sense of the order in the class of policies $\Pi$.

3.4. *Discussion.* The upper bounds established in Theorems 1 and 2 demonstrate that a finite regret can be achieved concurrent with an inferior sampling rate that grows to infinity. This is a rather obvious characteristic of the bandit problems with covariates: wrong arm "pulls" may incur a small, even negligible, loss in terms of rewards. In contrast, in traditional multi-armed bandit problems the regret and inferior sampling rate are identical up to a constant multiplier.

Woodroofe (1979) and Sarkar (1991) establish asymptotic optimality of the myopic policies in the Bayesian setting. In contrast, the rate-optimal policies developed here are nonmyopic; we were not able to show that the myopic policy is rate optimal in our setting. We believe the explanation for this lies in the following assumptions made in the aforementioned papers: Woodroofe (1979), Conditions C1 and C2, and Sarkar (1991), Conditions

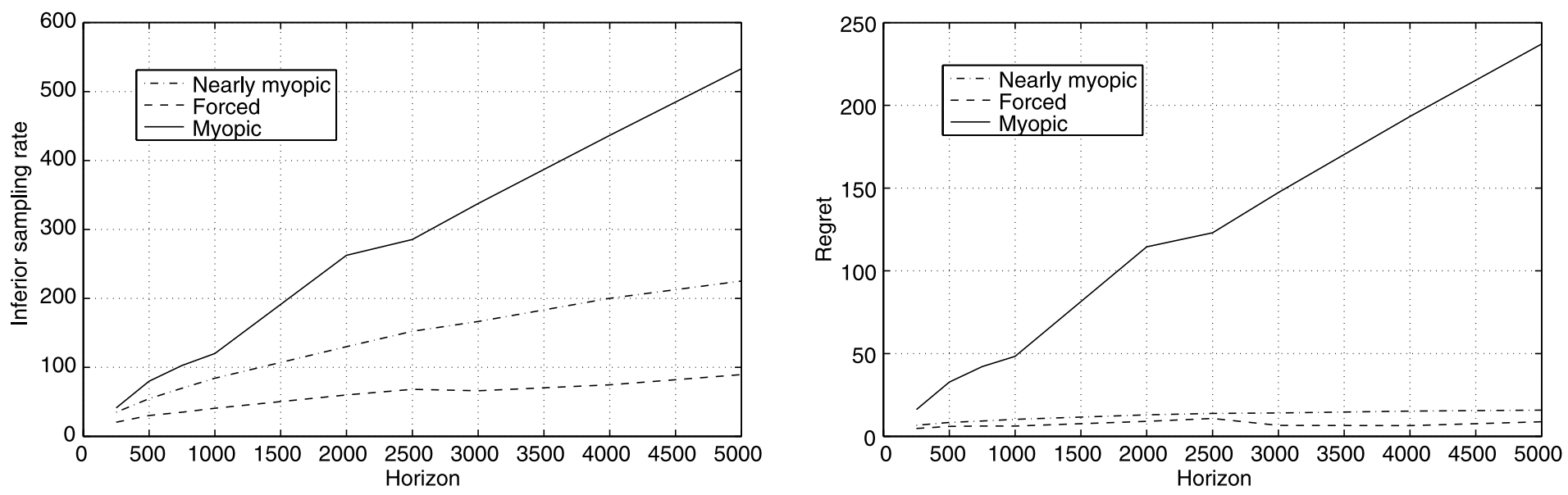

(a) (b)

FIG. 1. *Setup* (i)*: $X_t$ are i.i.d. uniformly distributed on* $[-1,1]$. (a) *The inferior sampling rate of the three policies averaged over 500 runs;* (b) *the regret of the three policies averaged over 500 runs.*

5a and 5b. These assumptions impose that the prior distribution for $\theta$ is supported on an interval $\Theta$, while the covariate $X$ has a positive continuous probability density on the real line. Therefore, with positive probability $X$ takes values outside $\Theta$, and for such covariate values it is exactly known which arm is superior (in expectation). This assumption ensures that for every $t$ with large probability the myopic rule samples $O(t)$ times from the arm 1 (cf. Lemmas 2 and 4). Note that the conditions (2.3) and (2.4) are much less restrictive. With the extra assumption made in Woodroofe (1979) and Sarkar (1991), one can establish optimality of a myopic policy, but it is worth pointing out that these assumptions simplify significantly the exploration–exploitation trade-off that underlies the design of "good" policies.

**4. Numerical results.** This section describes the results of a small simulation experiment that illustrates behavior of the policies presented in Section 3, and compares them with the myopic policy.

The conditional distribution of $Y$ given $X = x$ is assumed to be Gaussian with mean $x$ and variance 1; hence, $\theta = 0$. The following two setups are considered: (i) $X_t$ are i.i.d. random variables uniformly distributed on $[-1,1]$; (ii) $X_t$ are i.i.d. random variables taking values $\pm 1$ with probability 1/2. The former corresponds to a case where $\alpha = 1$, and the latter to a case where $\alpha = \infty$. In each setup we compute the inferior sampling rate and the regret of the three policies (myopic, nearly myopic and the one involving forced sampling), when the horizon $n$ takes values in the set $\{250, 500, 750, 1000, 2000, 2500, 3000, 4000, 5000\}$. In our simulations the nearly-myopic policy is implemented with $\delta_t = \sqrt{\ln t}$, while the forced sampling policy uses $q = 1/12$; see the conditions of Theorem 2. For each $n$ we compute the inferior sampling rate and the regret, averaged over 500 runs.

The results are summarized in Figures 1 and 2. Figure 1(a) shows the inferior sampling rate of the nearly myopic, forced sampling and myopic policies averaged over 500 runs, while Figure 1(b) displays the corresponding averaged regret. It is clearly seen that when the covariates $X_t$ are uniformly distributed, the forced sampling policy has the smallest averaged inferior sampling rate and regret. The nearly-myopic policy also outperforms the myopic policy which has the largest average inferior sampling rate and regret. Figure 2 corresponds to setup (ii) where $\alpha = \infty$. Because $X_t$ are i.i.d. random variables taking the values $\pm 1$ and $\theta = 0$, the inferior sampling rate coincides with the regret. That is why in Figure 2 we present only the graph of the logarithm of the average regret. The numerical results show that the nearly-myopic policy is preferable under setup (ii), consistent with the theoretical results of Section 3.

Even though the myopic policy appears to be inferior in comparison with the nearly-myopic and forced sampling policies, the results in Figures 1 and 2 do not clarify the reasons for such behavior. Additional insight into performance of the three policies can be gained from the graphs in Figure 3. For $n = 2000$ and under conditions of setup (i), Figures 3(a) and (b) display the boxplots of the inferior sampling rate and the regret obtained in 500 runs. It is clearly seen that the average performance of the myopic policy is badly affected by a large number of runs with poor performance. The nearly-myopic policy is the most stable over the different runs, though its average performance is worse than that of the forced sampling policy for this covariate distribution.

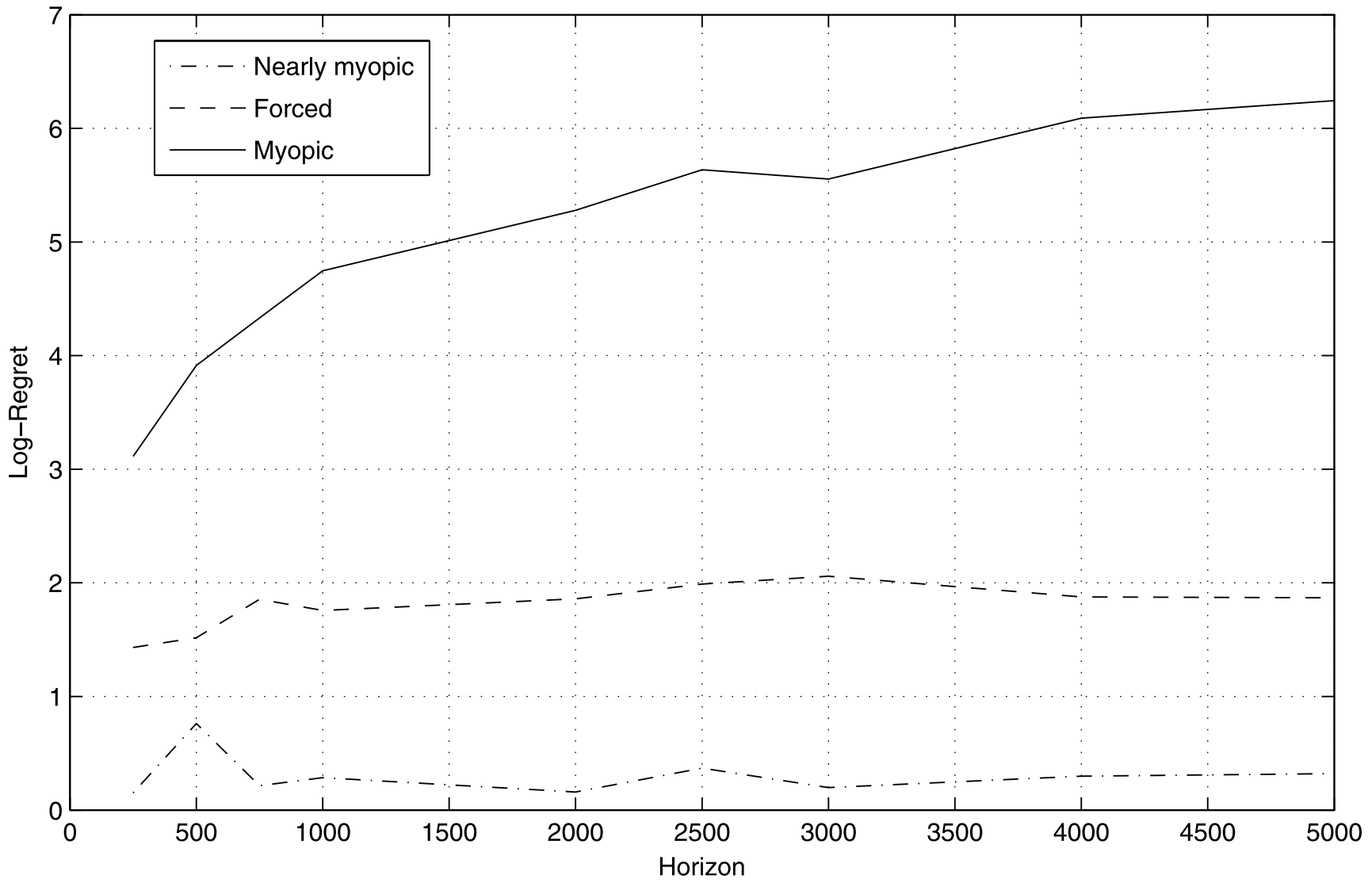

FIG. 2. *Setup* (ii)*:* $X_t$ *are i.i.d. random variables taking values* $\pm 1$ *with probability* $1/2$. *The logarithm of the regret averaged over 500 runs.*

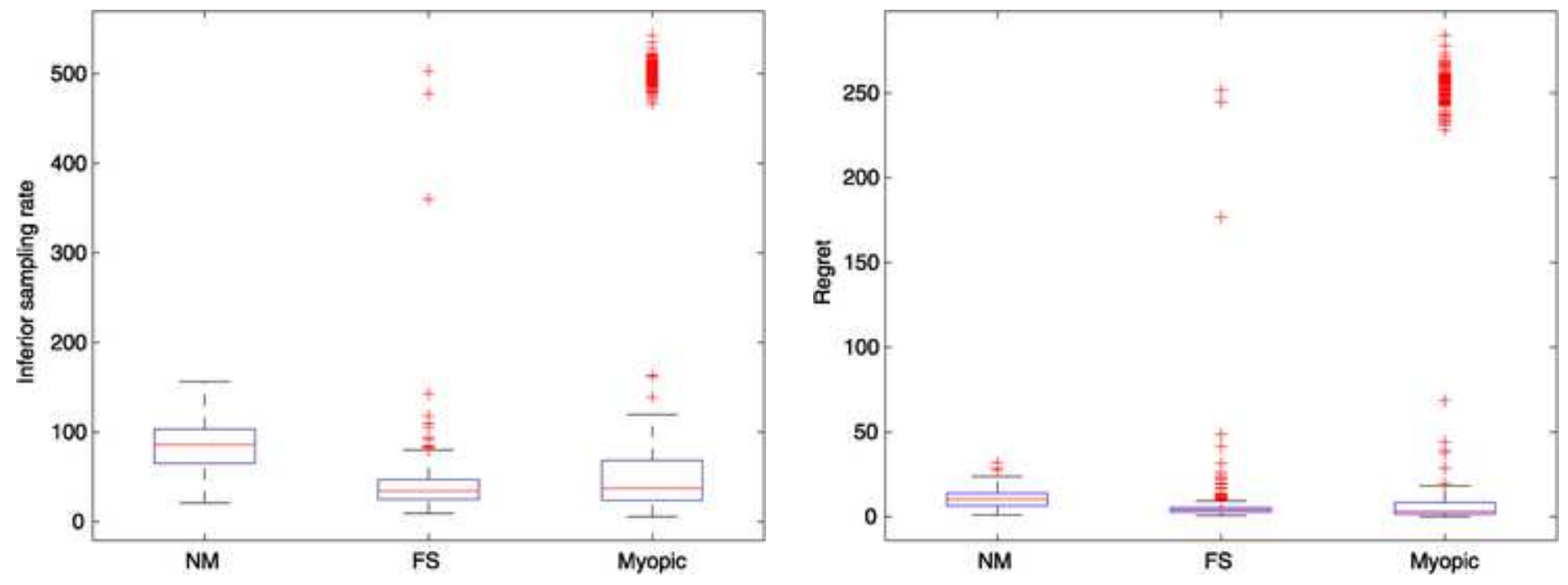

FIG. 3. *Setup* (i)*: $X_t$ are i.i.d. uniformly distributed on* $[-1,1]$. (a) *Boxplot of the inferior sampling rate computed over 500 runs;* (b) *Boxplot of the regret computed over 500 runs.*

## 5. Proofs.

### 5.1. *Preliminary lemma.*

LEMMA 1. *For any policy $\pi$, any measurable set $A$, and any $x>0$,*

$$\begin{aligned}&\mathbb{P}\{|\hat{\theta}_\pi(t)-\theta|>x,A\}\\&\quad\le 2\biggl[\mathbb{E}\exp\biggl\{-\frac{x^2}{2\sigma^2}T_\pi(t)\biggr\}I\{|\hat{\theta}_\pi(t)-\theta|>x,A\}\biggr]^{1/2}.\end{aligned}\tag{5.1}$$

*In particular,*

$$\mathbb{P}\{|\hat{\theta}_\pi(t)-\theta|>x,T_\pi(t)>\tau\}\le 2\exp\biggl\{-\frac{x^2\tau}{4\sigma^2}\biggr\}\qquad\forall x,\tau>0.\tag{5.2}$$

REMARK 5. The proof shows that (5.1) holds when $x$ is a positive random variable.

PROOF. The inequalities (5.1) and (5.2) follow immediately from results in de la Peña, Klass and Lai (2004) [see also de la Peña, Klass and Lai (2007) and Liptser and Spokoiny (2000)]. We provide a proof for completeness.

Note that $\hat{\theta}_\pi(t)-\theta=-(\sum_{s=1}^t\pi_s)^{-1}\sum_{s=1}^t\varepsilon_s\pi_s$. Write, for brevity, $S=\{\theta-\hat{\theta}_\pi(t)>x,A\}$. Then for any $\lambda>0$,

$$\mathbb{P}\{S\}=\mathbb{E}I\biggl\{\exp\biggl(\lambda\sum_{s=1}^t\varepsilon_s\pi_s-\lambda xT_\pi(t)\biggr)\ge 1\biggr\}I\{S\}$$

$$\leq \mathbb{E}\exp\left\{\lambda\sum_{s=1}^{t}\varepsilon_s\pi_s - \lambda x T_\pi(t)\right\} I\{S\}. \tag{5.3}$$

Let

$$M_t(\lambda) = \exp\left\{\lambda\sum_{s=1}^{t}\varepsilon_s\pi_s - \frac{\sigma^2\lambda^2}{2}T_\pi(t)\right\};$$

then $(M_t(\lambda), \mathcal{F}_t)$ is a martingale for any $\lambda$, and $\mathbb{E}M_t(\lambda) \leq 1$ for all $t$. Therefore, it follows from (5.3) that

$$\mathbb{P}\{S\} \leq \mathbb{E}\left(\exp\left\{\lambda\sum_{s=1}^{t}\varepsilon_s\pi_s - \lambda^2\sigma^2 T_\pi(t)\right\}\exp\{(\lambda^2\sigma^2 - \lambda x)T_\pi(t)\} I\{S\}\right)$$

$$\leq \sqrt{\mathbb{E}M_t(2\lambda)}(\mathbb{E}\exp\{2(\lambda^2\sigma^2 - \lambda x)T_\pi(t)\} I\{S\})^{1/2} \tag{5.4}$$

$$\leq \left[\mathbb{E}\exp\left\{-\frac{x^2}{2\sigma^2}T_\pi(t)\right\} I\{S\}\right]^{1/2},$$

where the second inequality is obtained using the Cauchy–Schwarz inequality, and the third one is by minimization over $\lambda > 0$. Applying the bound (5.4) for the random variable $-(\theta - \hat{\theta}_\pi(t))$, we complete the proof of the lemma. □

5.2. *Proof of Theorem 1.* We begin with the following lemma.

LEMMA 2. *Let (2.4) hold, and assume that the sequence* $\delta = (\delta_t, t \geq 1)$ *is nondecreasing. Then for any* $z \in (0, \frac{1}{4}]$ *and* $t \geq 2$,

$$\mathbb{P}\{T_{\hat{\pi}}(t) \leq zpt\} \leq \exp\left\{-\frac{1}{2}p^2z^2t\right\} + 4zt\exp\left\{-\frac{\delta^2_{\lfloor(1-2z)t\rfloor}}{4\sigma^2}\right\}. \tag{5.5}$$

PROOF. Denote $\zeta_s = T_{\hat{\pi}}^{-1}(s)\sum_{j=1}^{s}\varepsilon_j\hat{\pi}_j$. Then

$$\begin{aligned} T_{\hat{\pi}}(t) &= 1 + \sum_{s=2}^{t}\hat{\pi}_s \\ &= 1 + \sum_{s=2}^{t} I\{X_{s+1} \geq \hat{\theta}_{\hat{\pi}}(s) - \delta_s T_{\hat{\pi}}^{-1/2}(s)\} \\ &= 1 + \sum_{s=2}^{t} I\{X_{s+1} \geq \theta - \zeta_s - \delta_s T_{\hat{\pi}}^{-1/2}(s)\} \\ &\geq 1 + \sum_{s=2}^{t} I\{X_{s+1} \geq \theta\} I\{|\zeta_s| \leq \delta_s T_{\hat{\pi}}^{-1/2}(s)\}. \end{aligned}$$

Denote $w_s = I\{|\zeta_s| \le \delta_s T_{\hat\pi}^{-1/2}(s)\}$; then

$$\begin{aligned}
\mathbb{P}\{T_{\hat\pi}(t) \le zpt\} &\le \mathbb{P}\left\{\sum_{s=2}^t I(X_{s+1} \ge \theta) w_s \le zpt - 1\right\} \\
&\le \mathbb{P}\left\{\sum_{s=2}^t I(X_{s+1} \ge \theta) w_s \le zpt - 1, \sum_{s=2}^t w_s \ge 2zt\right\} \\
&\quad + \mathbb{P}\left\{\sum_{s=2}^t w_s < 2zt\right\} \\
&=: J_1 + J_2.
\end{aligned}$$

Let $p' = P_X\{[\theta, \infty)\}$; it follows from (2.4) that $p' \ge p$. Note that $w_s$ is $\mathcal{F}_{s-1}$-measurable; hence, $(\sum_{s=2}^t [p' - I\{X_{s+1} \ge \theta\}] w_s, \mathcal{F}_s)$ is the martingale with bounded differences. Then by the Azuma–Hoeffding inequality [see, e.g., Cesa-Bianchi and Lugosi (2006)],

$$\begin{aligned}
J_1 &= \mathbb{P}\left\{\sum_{s=2}^t [p' - I(X_{s+1} \ge \theta)] w_s \ge p' \sum_{s=2}^t w_s - (zpt - 1), \sum_{s=2}^t w_s \ge 2zt\right\} \\
&\le \mathbb{P}\left\{\sum_{s=2}^t [p' - I(X_{s+1} \ge \theta)] w_s \ge zpt + 1\right\} \\
&\le \exp\left\{-\frac{(zpt+1)^2}{2(t-2)}\right\} \le \exp\left\{-\frac{1}{2} z^2 p^2 t\right\}.
\end{aligned}$$

Now we bound $J_2$. For this purpose we note that

$$\left\{\sum_{s=2}^t I\{|\zeta_s| \ge \delta_s T_{\hat\pi}^{-1/2}(s)\} > (1-2z)t\right\} \subseteq \bigcup_{s=\lfloor (1-2z)t \rfloor}^t \{|\zeta_s| > \delta_s T_{\hat\pi}^{-1/2}(s)\}.$$

Therefore,

$$\begin{aligned}
J_2 &= \mathbb{P}\left\{\sum_{s=2}^t I\{|\zeta_s| > \delta_s T_{\hat\pi}^{-1/2}(s)\} > (1-2z)t\right\} \\
&\le \sum_{s=\lfloor (1-2z)t \rfloor}^t \mathbb{P}\{|\zeta_s| > \delta_s T_{\hat\pi}^{-1/2}(s)\} \\
&\le 2 \sum_{s=\lfloor (1-2z)t \rfloor}^t \exp\left\{-\frac{\delta_s^2}{4\sigma^2}\right\} \le 4zt \exp\left\{-\frac{\delta_{\lfloor (1-2z)t \rfloor}^2}{4\sigma^2}\right\},
\end{aligned}$$

where the second inequality follows form Lemma 1, and the third by monotonicity of $(\delta_t)$. Combining this inequality with the upper bound on $J_1$, we complete the proof. □

Now we turn to the proof of Theorem 1.

*Proof of Theorem 1.* (i) First we prove (3.4). Define

$$D_t = \left[\left(\hat{\theta}_{\hat{\pi}}(t) - \frac{\delta_t}{\sqrt{T_{\hat{\pi}}(t)}}\right) \wedge \theta, \left(\hat{\theta}_{\hat{\pi}}(t) - \frac{\delta_t}{\sqrt{T_{\hat{\pi}}(t)}}\right) \vee \theta\right].$$

Therefore,

$$T_{\inf}(n) = \sum_{t=1}^{n} I\{\hat{\pi}_t \neq \pi_t^*\} \leq 1 + \sum_{t=1}^{n-1} I\{X_{t+1} \in D_t\}$$

and

$$\begin{aligned} \mathbb{E}[T_{\inf}(n)] &\leq 1 + \sum_{t=1}^{n-1} \mathbb{P}\{X_{t+1} \in D_t\} \\ &\leq 1 + \sum_{t=1}^{n-1} [\mathbb{P}\{X_{t+1} \in D_t, T_{\hat{\pi}}(t) > pt/4\} + \mathbb{P}\{T_{\hat{\pi}}(t) \leq pt/4\}]. \end{aligned} \tag{5.6}$$

Applying Lemma 2 with $z = 1/4$, we have

$$\mathbb{P}\{T_{\hat{\pi}}(t) \leq pt/4\} \leq \exp\left\{-\frac{p^2 t}{32}\right\} + t \exp\left\{-\frac{\delta^2_{\lfloor t/2 \rfloor}}{4\sigma^2}\right\}. \tag{5.7}$$

Furthermore, for any sequence $(\gamma_t)$ of positive random variables such that $\gamma_t$ is $\mathcal{F}_{t-1}$-measurable, one can write

$$\begin{aligned} &\mathbb{P}\{X_{t+1} \in D_t, T_{\hat{\pi}}(t) > pt/4\} \\ &\quad = \mathbb{P}\{X_{t+1} \in D_t, T_{\hat{\pi}}(t) > pt/4, |\hat{\theta}_{\hat{\pi}}(t) - \delta_t T_{\hat{\pi}}^{-1/2}(t) - \theta| \leq \gamma_t\} \\ &\qquad + \mathbb{P}\{X_{t+1} \in D_t, T_{\hat{\pi}}(t) > pt/4, |\hat{\theta}_{\hat{\pi}}(t) - \delta_t T_{\hat{\pi}}^{-1/2}(t) - \theta| > \gamma_t\} \\ &\quad =: P_1(t) + P_2(t). \end{aligned}$$

Setting $\gamma_t = 2\delta_t T_{\hat{\pi}}^{-1/2}(t)$ and using the definition of $D_t$ and (2.3), we obtain

$$P_1(t) \leq \mathbb{P}\left\{|X_{t+1} - \theta| \leq \frac{2\delta_t}{\sqrt{T_{\hat{\pi}}(t)}}, T_{\hat{\pi}}(t) \geq pt/4\right\} \leq C_* \left(\frac{4\delta_t}{\sqrt{pt}}\right)^{\alpha}, \tag{5.8}$$

provided that $\delta_t \leq x_0 \sqrt{pt}/4$. By Lemma 1,

$$\begin{aligned} P_2(t) &\leq \mathbb{P}\{T_{\hat{\pi}}(t) > pt/4, |\hat{\theta}_{\hat{\pi}}(t) - \delta_t T_{\hat{\pi}}^{-1/2}(t) - \theta| > 2\delta_t T_{\hat{\pi}}^{-1/2}(t)\} \\ &\leq \mathbb{P}\{T_{\hat{\pi}}(t) > pt/4, |\hat{\theta}_{\hat{\pi}}(t) - \theta| > \delta_t T_{\hat{\pi}}^{-1/2}(t)\} \\ &\leq 2 \exp\left\{-\frac{\delta_t^2}{2\sigma^2}\right\}. \end{aligned} \tag{5.9}$$

Now we set $\delta_t = 2\sigma\sqrt{3\ln t}$; with this choice (5.9) and (5.7) imply that $P_2(t) \le 2t^{-6}$ and

$$\mathbb{P}\{T_{\hat\pi}(t) \le pt/4\} \le \exp\{-p^2 t/32\} + 8t^{-2} \qquad \forall t \ge 2.$$

Let $t_0 = t_0(p,\sigma,x_0) := \min\{t : x_0\sqrt{pt}/4 \ge 2\sigma\sqrt{3\ln t}\}$; then it follows from (5.8) that

$$P_1(t) \le C_*\left(8\sqrt{3}\sigma\sqrt{\frac{\ln t}{pt}}\right)^{\alpha} \qquad \forall t \ge t_0.$$

It is easily seen that $t_0 \le c_1\sigma^2(px_0^2)^{-1}\sqrt{|\ln[\sigma^2(px_0^2)^{-1}]|}$ for some absolute constant $c_1$. Combining these inequalities with (5.6), we obtain

$$\begin{aligned}\mathbb{E}[T_{\inf}(n)] \le (t_0 \vee 2) + C_* \sum_{t=t_0\vee 2}^{n-1}\left(8\sqrt{3}\sigma\sqrt{\frac{\ln t}{pt}}\right)^{\alpha} \\ + 8\sum_{t=t_0\vee 2}^{n-2} t^{-2} + 2\sum_{t=t_0\vee 2}^{n-1} t^{-6} + \sum_{t=t_0\vee 2}^{n-1}\exp\{-p^2t/32\}.\end{aligned}$$

This completes the proof of (3.4).

Now we prove (3.5). Assume that $\alpha > 2$ and let $t_\alpha$ be given by (3.3). Clearly,

$$t_\alpha \le c_2[(8\sqrt{3}\sigma p^{-1/2})\sqrt{\ln(8\sqrt{3}\sigma p^{-1/2})}]^{4\alpha/(\alpha-2)}\left[\frac{4\alpha}{\alpha-2}\right]^{2\alpha/(\alpha-2)}$$

for some absolute constant $c_2$. We can write

$$\mathbb{E}[T_{\inf}(n)] \le t_\alpha + \sum_{t=t_\alpha}^{n}\mathbb{P}\{\hat\pi_t \ne \pi_t^*\}.$$

Each summand on the right-hand side of the above formula is bounded from above exactly as before. This leads to the following bound:

$$\begin{aligned}\mathbb{E}[T_{\inf}(n)] \le (t_0 \vee 2 \vee t_\alpha) + C_* \sum_{t=t_0\vee 2\vee t_\alpha}^{n-1}\left(8\sqrt{3}\sigma\sqrt{\frac{\ln t}{pt}}\right)^{\alpha} \\ + 8\sum_{t=t_0\vee 2\vee t_\alpha}^{n-2} t^{-2} + 2\sum_{t=t_0\vee 2\vee t_\alpha}^{n-1} t^{-6} \\ + \sum_{t=t_0\vee 2\vee t_\alpha}^{n-1}\exp\{-p^2t/32\}.\end{aligned} \tag{5.10}$$

By definition of $t_\alpha$, $(8\sqrt{3}\sigma\sqrt{\ln t/p})^\alpha \le t^{\alpha/4-1/2}$ for all $t \ge t_\alpha$; hence,

$$\sum_{t=t_\alpha}^{n-1}\left(8\sqrt{3}\sigma\sqrt{\frac{\ln t}{pt}}\right)^\alpha = \sum_{t=t_\alpha}^{n-1}\left(\frac{8\sqrt{3}\sigma}{\sqrt{p}}\right)^\alpha \frac{(\ln t)^{\alpha/2}}{t^{\alpha/4-1/2}}\frac{1}{t^{\alpha/4+1/2}}$$
$$\le \sum_{t=t_\alpha}^{n-1} t^{-(\alpha+2)/4} \le \frac{4n^{-4/(\alpha-2)}}{\alpha-2}.$$

This bound along with (5.10) leads to (3.5).

(ii) The proof of the second statement goes along the same lines.

We have

$$\begin{aligned}
&\sum_{t=1}^{n}|X_t-\theta| I\{\hat\pi_t \ne \pi_t^*\}\\
&\quad\le |X_1-\theta| + \sum_{t=1}^{n-1}|X_{t+1}-\theta| I\{X_{t+1}\in D_t\}\\
&\quad\le |X_1-\theta| + \sum_{t=1}^{n-1}|X_{t+1}-\theta| I\{X_{t+1}\in D_t, T_{\hat\pi}(t) > pt/4\}\\
&\qquad + \sum_{t=1}^{n-1}|X_{t+1}-\theta| I\{T_{\hat\pi}(t) \le pt/4\}\\
&\quad =: |X_1-\theta| + J_1 + J_2.
\end{aligned}$$

Since $\mathbb{E}|X_t-\theta| \le \mu$ for $P_{X,Y} \in \mathcal{P}'_\alpha$, and $X_{t+1}$ is independent of $T_{\hat\pi}(t)$, we have by (5.7)

$$\mathbb{E}[J_2] \le \mu\sum_{t=1}^{n-1}\left[\exp\left\{-\frac{p^2 t}{32}\right\} + t\exp\left\{-\frac{\delta^2_{\lfloor t/2\rfloor}}{4\sigma^2}\right\}\right].$$

Furthermore,

$$\begin{aligned}
\mathbb{E}[J_1] &= \sum_{t=1}^{n-1}\mathbb{E}|X_{t+1}-\theta| I\{X_{t+1}\in D_t, T_{\hat\pi}(t) > pt/4,\\
&\qquad\qquad |\hat\theta_{\hat\pi}(t) - \delta_t T_{\hat\pi}^{-1/2}(t) - \theta| \le \gamma_t\}\\
&\quad + \sum_{t=1}^{n-1}\mathbb{E}|X_{t+1}-\theta| I\{X_{t+1}\in D_t, T_{\hat\pi}(t) > pt/4,\\
&\qquad\qquad |\hat\theta_{\hat\pi}(t) - \delta_t T_{\hat\pi}^{-1/2}(t) - \theta| > \gamma_t\}\\
&=: \sum_{t=1}^{n-1}E_1(t) + \sum_{t=1}^{n-1}E_2(t).
\end{aligned}$$

Setting as before $\gamma_t = 2\delta_t T_{\hat{\pi}}^{-1/2}(t)$ and using (5.8), we obtain

$$E_1(t) \le \mathbb{E}|X_{t+1}-\theta| I\left\{|X_{t+1}-\theta| \le \frac{2\delta_t}{\sqrt{T_{\hat{\pi}}(t)}}, T(t) \ge pt/4\right\} \le C_*\left(\frac{4\delta_t}{\sqrt{pt}}\right)^{\alpha+1},$$

provided that $\delta_t \le x_0\sqrt{pt}/4$. In addition,

$$\begin{aligned} E_2(t) &= \mathbb{E}|X_{t+1}-\theta| I\{X_{t+1} \in D_t, T_{\hat{\pi}}(t) > pt/4, |\hat{\theta}_{\hat{\pi}}(t) - \delta_t T_{\hat{\pi}}^{-1/2}(t) - \theta| > \gamma_t\} \\ &\le \mathbb{E}|X_{t+1}-\theta| I\{T_{\hat{\pi}}(t) > pt/4, |\hat{\theta}_{\hat{\pi}}(t) - \delta_t T_{\hat{\pi}}^{-1/2}(t) - \theta| > \gamma_t\} \\ &\le 2\mu \exp\left\{-\frac{\delta_t^2}{2\sigma^2}\right\}. \end{aligned}$$

Combining these inequalities, we come to (3.6). The bound (3.7) is obtained using the same reasoning as in the proof of (3.5).

5.3. *Proof of Theorem 2.* The next result follows from Lemma 1 by setting $A = \Omega$ and taking into account that $T_{\tilde{\pi}}(t) \ge N(t)$, $\forall t$.

LEMMA 3. *Under policy $\tilde{\pi}$ for any $x > 0$ and any $t \ge 1$, one has*

$$\mathbb{P}\{|\hat{\theta}_{\tilde{\pi}}(t) - \theta| > x\} \le 2\exp\left\{-\frac{x^2}{4\sigma^2} N(t)\right\}.$$

In particular, it follows from (3.12) and Lemma 3 that

$$\mathbb{P}\{|\hat{\theta}_{\tilde{\pi}}(t) - \theta| > x\} \le 2\exp\left\{\frac{x^2}{4\sigma^2}\right\}\left(\frac{t}{\nu+1}\right)^{-x^2/(4q\sigma^2)} \qquad \forall t > \nu, \tag{5.11}$$

where $\nu$ is given in (3.10).

LEMMA 4. *Let $P_{X,Y} \in \mathcal{P}_\alpha$ and $p_1 := p - C_* x_0^\alpha$; then for any $z \in (0, 1/4]$ and all $t > \nu$,*

$$\begin{aligned} &\mathbb{P}\{T_{\tilde{\pi}}(t) \le z p_1 t\} \\ &\quad \le \exp\left\{-\frac{1}{2} p_1^2 z^2 (t - N(t))\right\} \\ &\qquad + 2\exp\left\{\frac{x_0^2}{4\sigma^2}[1 + q^{-1}\ln(\nu+1)]\right\} t^{1 - x_0^2/(4q\sigma^2)}. \end{aligned}$$

PROOF. Fix $t > \nu$. Denote $\zeta_t = \frac{1}{T_{\tilde{\pi}}(t)} \sum_{j=1}^t \varepsilon_j \tilde{\pi}_j$. Then

$$T_{\tilde{\pi}}(t) = N(t) + \sum_{s=1}^t \tilde{\pi}_s I\{s \in \mathcal{T}^c\}$$

$$
\begin{aligned}
&= N(t) + \sum_{s=1}^{t} I\{X_{s+1} \geq \hat{\theta}_{\tilde{\pi}}(s)\} I\{s \in \mathcal{T}^c\} \\
&= N(t) + \sum_{s=1}^{t} I\{X_{s+1} \geq \theta - \zeta_s\} I\{s \in \mathcal{T}^c\} \\
&\geq N(t) + \sum_{s=1}^{t} I\{X_{s+1} \geq \theta + x_0\} I\{|\zeta_s| \leq x_0\} I\{s \in \mathcal{T}^c\}.
\end{aligned}
$$

Denote $w_s = I\{|\zeta_s| \leq x_0\} I\{s \in \mathcal{T}^c\}$. Then

$$
\begin{aligned}
&\mathbb{P}\{T_{\tilde{\pi}}(t) \leq z p_1 t\} \\
&\quad \leq \mathbb{P}\left\{\sum_{s=1}^{t} I\{X_{s+1} \geq \theta + x_0\} w_s \leq z p_1 (t - N(t))\right\} \\
&\quad \leq \mathbb{P}\left\{\sum_{s=1}^{t} I\{X_{s+1} \geq \theta + x_0\} w_s \leq z p_1 (t - N(t)), \sum_{s=1}^{t} w_s \geq 2z(t - N(t))\right\} \\
&\qquad + \mathbb{P}\left\{\sum_{s=1}^{t} w_s < 2z(t - N(t))\right\} \\
&\quad =: J_1(t) + J_2(t).
\end{aligned}
$$

Let $p_1' = \mathbb{P}\{X_{s+1} \geq \theta + x_0\}$; it follows from the definition of $\mathcal{P}_\alpha$ [cf. (2.5)] that $p_1' \geq p_1 > 0$. Note that $w_s$ is $\mathcal{F}_{s-1}$-measurable, and $(\sum_{s=1}^{t} [p_1' - I\{X_{s+1} \geq \theta + x_0\}] w_s, \mathcal{F}_s)$ is a martingale with bounded differences. Then by the Azuma–Hoeffding inequality [see, e.g., Cesa-Bianchi and Lugosi (2006)],

$$
\begin{aligned}
J_1(t) &= \mathbb{P}\left\{\sum_{s=1}^{t} [p_1' - I\{X_{s+1} \geq \theta + x_0\}] w_s \geq p_1' \sum_{s=1}^{t} w_s - z p_1' (t - N(t)),\right. \\
&\qquad\qquad \left. \sum_{s=1}^{t} w_s \geq 2z(t - N(t))\right\} \\
&\leq \mathbb{P}\left\{\sum_{s=1}^{t} [p_1' - I\{X_{s+1} \geq \theta + x_0\}] w_s \geq z p_1 (t - N(t))\right\} \\
&\leq \exp\left\{-\frac{1}{2} z^2 p_1^2 (t - N(t))\right\}.
\end{aligned}
$$

Now we bound $J_2(t)$ as follows:

$$
J_2(t) = \mathbb{P}\left\{\sum_{s=1}^{t} w_s < 2z(t - N(t))\right\}
$$

$$
\begin{aligned}
&= \mathbb{P}\Bigg\{\sum_{s=1}^{t} I\{|\zeta_s| \le x_0\} I\{s \in \mathcal{T}^c\} < 2z(t - N(t))\Bigg\} \\
&= \mathbb{P}\Bigg\{\sum_{s=1}^{t} I\{|\zeta_s| > x_0\} I\{s \in \mathcal{T}^c\} \ge (1-2z)(t - N(t))\Bigg\} \\
&\le \mathbb{P}\Bigg\{\bigcup_{s=\lfloor (t-N(t))(1-2z)\rfloor}^{t} \{|\zeta_s| > x_0\}\Bigg\} \\
&\le (t - N(t)) \max_{t-N(t)\le s \le t} \mathbb{P}\{|\zeta_s| \ge x_0\} \\
&\le 2\exp\Bigg\{\frac{x_0^2}{4\sigma^2}[1 + q^{-1}\ln(\nu+1)]\Bigg\} t^{1-x_0^2/(4q\sigma^2)},
\end{aligned}
$$

where the last inequality follows from (5.11). Combining this inequality with the upper bound on $J_1(t)$ completes the proof. □

*Proof of Theorem 2.* $1^0$. First we prove (3.15). By the premise of the theorem,

$$x_0^2/(12\sigma^2) \ge q. \tag{5.12}$$

Put $\bar{\mathcal{T}}_t = \{1 \le s < t : s \notin \mathcal{T}\}$, and define $D_t = [\hat{\theta}_{\tilde{\pi}}(t) \wedge \theta, \hat{\theta}_{\tilde{\pi}}(t) \vee \theta]$. Then

$$T_{\inf}(n) = \sum_{t=1}^{n} I\{\tilde{\pi}_t \ne \pi_t^*\} \le N(n) + \sum_{t \in \bar{\mathcal{T}}_n} I\{X_{t+1} \in D_t\}.$$

By choice of the "forced" sampling sequence $\mathcal{T}$ in view of (3.11), we have that

$$
\begin{aligned}
\mathbb{E}[T_{\inf}(n)] &\le 1 + \frac{1}{q}\ln(n+1) + \sum_{t\in\bar{\mathcal{T}}_n} \mathbb{P}\{X_{t+1} \in D_t\} \\
&= 1 + \frac{1}{q}\ln(n+1) + \sum_{t\in\bar{\mathcal{T}}_n} \mathbb{P}\{X_{t+1} \in D_t, T_{\tilde{\pi}}(t) > p_1 t/4\} \\
&\quad + \sum_{t\in\bar{\mathcal{T}}_n} \mathbb{P}\{X_{t+1} \in D_t, T_{\tilde{\pi}}(t) \le p_1 t/4\} \\
&\le 1 + \frac{1}{q}\ln(n+1) + \nu_0 + \sum_{t\in\bar{\mathcal{T}}_n} P_1(t) + \sum_{t\in\bar{\mathcal{T}}_n, t>\nu_0} P_2(t),
\end{aligned} \tag{5.13}
$$

where $\nu_0$ is defined in (3.14). Applying Lemma 4 with $z = 1/4$, we have

$$P_2(t) \le \mathbb{P}\{T_{\tilde{\pi}}(t) \le p_1 t/4\}$$

$$\le \exp\left\{-\frac{p_1^2}{32}(t-N(t))\right\} + 2\exp\left\{\frac{x_0^2}{4\sigma^2}[1+q^{-1}\ln(\nu+1)]\right\}t^{1-x_0^2/(4q\sigma^2)}.$$

$$\le \exp\{-p_1^2 t/64\} + 2\exp\left\{\frac{x_0^2}{4\sigma^2}[1+q^{-1}\ln(\nu+1)]\right\}t^{-2},$$

where the last inequality follows from (5.12) and the fact that $t-N(t)\ge t/2$ for $t>\nu_0$. Hence, we get that

$$\begin{aligned}\sum_{t\in\bar{\mathcal{T}}_n, t>\nu_0} P_2(t) \le{}& \frac{1}{1-\exp\{-p_1^2/64\}}\\ &+\frac{\pi^2}{3}\exp\left\{\frac{x_0^2}{4\sigma^2}[1+q^{-1}\ln(\nu+1)]\right\}.\end{aligned}\tag{5.14}$$

Now, turning to $P_1(t)$, we have

$$\begin{aligned}P_1(t) &= \mathbb{P}\{X_{t+1}\in D_t, T_{\tilde\pi}(t)>p_1t/4\}\\ &= \mathbb{P}\{X_{t+1}\in D_t, T_{\tilde\pi}(t)>p_1t/4, |\hat\theta_{\tilde\pi}(t)-\theta|\le x_0\}\\ &\quad+\mathbb{P}\{X_{t+1}\in D_t, T_{\tilde\pi}(t)>p_1t/4, |\hat\theta_{\tilde\pi}(t)-\theta|> x_0\}\\ &=: J_1(t)+J_2(t).\end{aligned}$$

We first bound $J_2(t)$. Using Lemma 1, we have

$$J_2(t)\le \mathbb{P}\{T_{\tilde\pi}(t)>p_1t/4, |\hat\theta_{\tilde\pi}(t)-\theta|>x_0\}\le 2\exp\left\{-\frac{x_0^2p_1t}{8\sigma^2}\right\}$$

and, therefore,

$$\sum_{t\in\bar{\mathcal{T}}_n} J_2(t)\le\sum_{t=1}^{n}2\exp\left\{-\frac{x_0^2p_1t}{8\sigma^2}\right\}=\frac{2}{1-\exp\{-x_0^2p_1/(8\sigma^2)\}}.\tag{5.15}$$

For $J_1(t)$, we proceed as follows:

$$\begin{aligned}J_1(t) &= \sum_{k=0}^{\infty}\mathbb{P}\{X_{t+1}\in D_t, T_{\tilde\pi}(t)>p_1t/4, 2^{-k-1}x_0<|\hat\theta_{\tilde\pi}(t)-\theta|\le 2^{-k}x_0\}\\ &= \sum_{k=0}^{\infty}\mathbb{E}[I\{T_{\tilde\pi}(t)>p_1t/4, 2^{-k-1}x_0<|\hat\theta_{\tilde\pi}(t)-\theta|\le 2^{-k}x_0\}\\ &\qquad\qquad\times\mathbb{P}\{X_{t+1}\in D_t|\mathcal{F}_t\}]\\ &\overset{(a)}{\le}\sum_{k=0}^{\infty}C_*[2^{-k}x_0]^{\alpha}\mathbb{P}\{T_{\tilde\pi}(t)>p_1t/4, |\hat\theta_{\tilde\pi}(t)-\theta|>2^{-k-1}x_0\}\\ &\overset{(b)}{\le}\sum_{k=0}^{\infty}C_*[2^{-k}x_0]^{\alpha}2\exp\left\{-\frac{x_0^22^{-2k-2}}{2\sigma^2}\frac{p_1t}{4}\right\}\end{aligned}\tag{5.16}$$

$$= 2C_* x_0^\alpha \sum_{k=0}^{\infty} 2^{-\alpha k} \exp\left\{-\frac{2^{-2k} x_0^2 t p_1}{32\sigma^2}\right\},$$

where (a) follows from condition (2.3) and (b) follows from Lemma 1.

For $\alpha, b > 0$, set

$$\begin{aligned} S(\alpha, b) &:= \sum_{k=0}^{\infty} 2^{-\alpha k} \exp\{-b 2^{-2k}\}, \\ I(\alpha, b) &:= \int_0^\infty 2^{-\alpha y} \exp\{-b 2^{-2y}\}\, dy. \end{aligned} \tag{5.17}$$

Note that the integrand above has a unique (global) maximum at $y^* = -(1/2) \times \log_2(\alpha/2b)$, provided that $\alpha \le 2b$. Put $k^* := \lfloor y^* \rfloor$ and write

$$\begin{aligned} S(\alpha, b) &= \sum_{k=0}^{k^*} 2^{-\alpha k} \exp\{-b 2^{-2k}\} + \sum_{k=k^*+1}^{\infty} 2^{-\alpha k} \exp\{-b 2^{-2k}\} \\ &=: S_1(\alpha, b) + S_2(\alpha, b). \end{aligned}$$

It follows that

$$S_2(\alpha, b) \le \frac{2^{-\alpha y^*}}{1 - 2^{-\alpha}} = \frac{(\alpha/2)^{\alpha/2}}{1 - 2^{-\alpha}} b^{-\alpha/2}.$$

Since the integrand in (5.17) is monotone increasing on $[0, y^*)$, we have that

$$\begin{aligned} S_1(\alpha, b) &\le \int_0^{y^*} 2^{-\alpha y} \exp\{-b 2^{-2y}\}\, dy \le I(\alpha, b) \\ &= \frac{1}{\ln 2} \int_0^1 z^{\alpha - 1} \exp\{-b z^2\}\, dz \le \frac{\Gamma(\alpha/2)}{2 \ln 2} b^{-\alpha/2}, \end{aligned}$$

where $\Gamma(\cdot)$ denotes the gamma function. Thus, we have shown that for all $0 < \alpha \le 2b$ one has

$$S(\alpha, b) \le b^{-\alpha/2} \left[\frac{(\alpha/2)^{\alpha/2}}{1 - 2^{-\alpha}} + \frac{\Gamma(\alpha/2)}{2 \ln 2}\right].$$

Now we apply this result with $b = x_0^2 t p_1/(32\sigma^2)$ in order to bound $J_2(t)$ [see (5.16)]. In particular, for any $t \ge 16\sigma^2 \alpha/(x_0^2 p_1)$, we have

$$J_1(t) \le 6C_* \left(\frac{32\sigma^2}{t p_1}\right)^{\alpha/2} \left[\frac{(\alpha/2)^{\alpha/2}}{1 - 2^{-\alpha}} + \frac{\Gamma(\alpha/2)}{2 \ln 2}\right],$$

and, hence,

$$\sum_{t \in \bar{\mathcal{T}}_n} J_1(t) \le \frac{16 \alpha \sigma^2}{x_0^2 p_1} + 6C_* \left[\frac{(\alpha/2)^{\alpha/2}}{1 - 2^{-\alpha}} + \frac{\Gamma(\alpha/2)}{2 \ln 2}\right] \sum_{t=1}^{n} \left(\frac{32\sigma^2}{t p_1}\right)^{\alpha/2}. \tag{5.18}$$

Combining (5.18) with (5.15), (5.14) and (5.13), we come to (3.15).

Now consider the case of $\alpha > 2$. Here we have

$$\sum_{t\in\bar{\mathcal{T}}_n} J_1(t) \le 2C_* x_0^\alpha \sum_{t=1}^{n}\sum_{k=0}^{\infty} 2^{-\alpha k}\exp\biggl\{-\frac{2^{-2k}x_0^2 t p_1}{32\sigma^2}\biggr\}$$
$$\le 2C_* x_0^\alpha \sum_{k=0}^{\infty}\frac{2^{-\alpha k}}{1-\exp\{-2^{-2k}x_0^2 p_1/(32\sigma^2)\}}.$$

If $x_0^2 p_1/(32\sigma^2) > 1$, then for all $k > k_0 = (2\ln 2)^{-1}\ln(x_0^2 p_1/32\sigma^2)$ we have that $2^{-2k}x_0^2 p_1/(32\sigma^2) \le 1$. The last inequality holds for all $k$ if $x_0^2 p_1/(32\sigma^2) \le 1$. In both cases

$$\sum_{k=0}^{\infty}\frac{2^{-\alpha k}}{1-\exp\{-2^{-2k}x_0^2 p_1/(32\sigma^2)\}}$$
$$\le \frac{1}{1-e^{-1}}\sum_{k=0}^{\lfloor k_0\rfloor} 2^{-\alpha k}$$
$$+ \sum_{k=\lfloor k_0\rfloor+1}^{\infty}\frac{2^{-\alpha k}}{2^{-2k}x_0^2 p_1/(32\sigma^2) - 1/2[2^{-2k}x_0^2 p_1/(32\sigma^2)]^2}$$
$$\le \frac{1}{(1-e^{-1})(1-2^{-\alpha})} + \frac{64\sigma^2}{x_0^2 p_1}\sum_{k=\lfloor k_0\rfloor+1}^{\infty} 2^{-(\alpha-2)k}$$
$$\le \frac{1}{(1-e^{-1})(1-2^{-\alpha})} + \frac{2}{1-2^{2-\alpha}}\biggl(\frac{32\sigma^2}{x_0^2 p_1}\biggr)^{\alpha/(\alpha-2)},$$

where in the second inequality we took into account that $2^{-2k}x_0^2 p_1/(32\sigma^2) \le 1$ for $k > k_0$. Therefore, if $\alpha > 2$, then

$$\sum_{t\in\bar{\mathcal{T}}_n} J_1(t) \le 2C_* x_0^\alpha\biggl\{\frac{2}{1-2^{2-\alpha}}\biggl(\frac{32\sigma^2}{x_0^2 p_1}\biggr)^{\alpha/(\alpha-2)} + \frac{4}{3(1-e^{-1})}\biggr\}. \tag{5.19}$$

Combining (5.19) with (5.15), (5.14) and (5.13), we come to (3.16).

$2^0$. The proof of the second statement proceeds using almost identical arguments. We have

$$\sum_{t\in\bar{\mathcal{T}}_n} |X_t - \theta| I\{\tilde{\pi}_t \ne \pi_t^*\}$$
$$\le \sum_{t\in\bar{\mathcal{T}}_n} |X_{t+1} - \theta| I\{X_{t+1}\in D_t\} + \sum_{t\in\bar{\mathcal{T}}_n^c} |X_t - \theta|$$
$$\le \sum_{t\in\bar{\mathcal{T}}_n} |X_{t+1} - \theta| I\{X_{t+1}\in D_t, T_{\tilde{\pi}}(t) > p_1 t/4\}$$

$$
\begin{aligned}
&+ \sum_{t\in\bar{\mathcal{T}}_n, t>\nu_0} |X_{t+1}-\theta| I\{T_{\tilde{\pi}}(t) \le p_1 t/4\} \\
&+ \biggl( \sum_{t\in\bar{\mathcal{T}}_n, t\le\nu_0} |X_t-\theta| + \sum_{t\in\bar{\mathcal{T}}_n^c} |X_t-\theta| \biggr) \\
&=: J_1(n) + J_2(n) + J_3(n).
\end{aligned}
$$

Because $P_{X,Y} \in \mathcal{P}'_\alpha(\theta)$, $\mathbb{E}|X_t-\theta| \le \mu$. Then, using the properties of the "forced" sampling sequence $\mathcal{T}$, we have that

$$
\mathbb{E}[J_3(n)] \le \mu \biggl[ 1 + \nu_0 + \frac{1}{q}\ln(n+1) \biggr].
$$

Now, since $X_{t+1}$ is independent of $T_{\tilde{\pi}}(t)$, arguing in the same way as in (5.14), we have

$$
\mathbb{E}[J_2(n)] \le \mu \biggl[ \frac{1}{1-\exp\{-p_1^2/64\}} + \frac{\pi^2}{3}\exp\biggl\{ \frac{x_0^2}{4\sigma^2}[1+q^{-1}\ln(\nu+1)] \biggr\} \biggr].
$$

Furthermore,

$$
\begin{aligned}
\mathbb{E}[J_1(n)] &= \sum_{t\in\bar{\mathcal{T}}_n} \mathbb{E}|X_{t+1}-\theta| I\{X_{t+1}\in D_t, T_{\tilde{\pi}}(t) > p_1 t/4, |\hat{\theta}_{\tilde{\pi}}(t)-\theta| \le x_0\} \\
&\quad + \sum_{t\in\bar{\mathcal{T}}_n} \mathbb{E}|X_{t+1}-\theta| I\{X_{t+1}\in D_t, T_{\tilde{\pi}}(t) > p_1 t/4, |\hat{\theta}_{\tilde{\pi}}(t)-\theta| > x_0\} \\
&=: \sum_{t=1}^{n-1} E_1(t) + \sum_{t=1}^{n-1} E_2(t).
\end{aligned}
$$

Using Lemma 1 and the independence of $X_{t+1}$ from $T_{\tilde{\pi}}(t), \hat{\theta}_{\tilde{\pi}}(t)$, we bound the second term as follows:

$$
\begin{aligned}
E_2(t) &= \mathbb{E}|X_{t+1}-\theta| I\{X_{t+1}\in D_t, T_{\tilde{\pi}}(t) > p_1 t/4, |\hat{\theta}_{\tilde{\pi}}(t)-\theta| > x_0\} \\
&\le \mathbb{E}|X_{t+1}-\theta| I\{T_{\tilde{\pi}}(t) > p_1 t/4, |\hat{\theta}_{\tilde{\pi}}(t)-\theta| > x_0\} \\
&\le 2\mu \exp\biggl\{ -\frac{x_0^2 p_1 t}{8\sigma^2} \biggr\}.
\end{aligned}
$$

Now, for the first term, write

$$
\begin{aligned}
E_1(t) &= \sum_{k=0}^{\infty} \mathbb{E}[|X_{t+1}-\theta| \\
&\qquad \times I\{X_{t+1}\in D_t, T_{\tilde{\pi}}(t) > p_1 t/4, 2^{-k-1}x_0 < |\hat{\theta}_{\tilde{\pi}}(t)-\theta| \le 2^{-k}x_0\}] \\
&= \sum_{k=0}^{\infty} \mathbb{E}[I\{T_{\tilde{\pi}}(t) > p_1 t/4, 2^{-k-1}x_0 < |\hat{\theta}_{\tilde{\pi}}(t)-\theta| \le 2^{-k}x_0\}
\end{aligned}
$$

$$\begin{aligned}
&\times \mathbb{E}\{|X_{t+1}-\theta| I\{X_{t+1}\in D_t\}|\mathcal{F}_t\}] \\
&\overset{\text{(a)}}{\le} \sum_{k=0}^{\infty} C_*[2^{-k}x_0]^{\alpha+1}\mathbb{P}\{T_{\tilde{\pi}}(t) > p_1 t/4, |\hat{\theta}_{\tilde{\pi}}(t)-\theta| > 2^{-k-1}x_0\} \\
&\overset{\text{(b)}}{\le} 2C_* x_0^{\alpha+1}\sum_{k=0}^{\infty} 2^{-(\alpha+1)k}\exp\left\{-\frac{2^{-2k}x_0^2 t p_1}{32\sigma^2}\right\} \\
&= 2C_* x_0^{\alpha+1} S(\alpha+1, b),
\end{aligned}$$

where, (a) follows from condition (2.3), (b) follows from Lemma 1, and with the notation used earlier, in part $1^0$, the function $S(\cdot,\cdot)$ is defined in (5.17) and $b := x_0^2 p_1 t/(32\sigma^2)$. Using the bound on $S(\alpha, b)$ derived earlier, and making the substitution $\alpha \mapsto (\alpha+1)$, we get

$$E_1(t) \le 2C_*\left[\frac{((\alpha+1)/2)^{(\alpha+1)/2}}{1-2^{-(\alpha+1)/2}} + \frac{\Gamma((\alpha+1)/2)}{2\ln 2}\right]\left(\frac{32\sigma^2}{tp_1}\right)^{(\alpha+1)/2}.$$

Summing over $t$ and using the bounds derived above on $E_2(t)$ together with the bounds established already on $\mathbb{E}[J_1(n)]$ and $\mathbb{E}[J_2(n)]$, we obtain (3.17).

If $\alpha > 1$, then arguing as in the proof of (3.15) we arrive at the result stated in (3.18).

5.4. *Proof of Theorem 3.* The proof relies on the following lemma.

LEMMA 5. *Let (2.3) hold; then for any policy $\pi$, one has*

$$R_n(\pi, \pi^*) \ge \frac{[S_n(\pi, \pi^*)]^{1+1/\alpha} n^{-1/\alpha}}{2\max\{(1/x_0), (2C_*)^{1/\alpha}\}}.$$

PROOF. Write, for brevity, $d_t(\pi, \pi^*) = \mathbb{P}\{\pi_t \ne \pi_t^*\}$. In order to underline dependence of $\pi_t$ on the observations $\mathcal{Y}_{t-1} = (X_1, \pi_1, \pi_1 Y_1, \ldots, X_{t-1}, \pi_{t-1}, \pi_{t-1} \times Y_{t-1})$ and on the covariate value $X_t$, we will write $\pi_t = \pi_t(\mathcal{Y}_{t-1}; X_t)$. We write also $\pi_t^* = \pi_t^*(X_t)$.

Let $\eta_t$ be a sequence of positive real numbers such that $\eta_t \le x_0$, $\forall t$; then

$$\begin{aligned}
R_n(\pi, \pi^*) &\ge \sum_{t=1}^{n} \mathbb{E}|X_t - \theta| I\{\pi_t(\mathcal{Y}_{t-1}; X_t) \ne \pi_t^*(X_t)\} I\{|X_t-\theta| > \eta_t\} \\
&\ge \sum_{t=1}^{n} \eta_t \int_{\{x : |x-\theta| > \eta_t\}} \mathbb{P}\{\pi_t(\mathcal{Y}_{t-1}; x) \ne \pi_t^*(x) | X_t = x\} P_X(dx) \\
&= \sum_{t=1}^{n} \eta_t \Big[ d_t(\pi, \pi^*)
\end{aligned}$$

(5.20)

$$- \int_{\{x\,:\,|x-\theta|\le\eta_t\}} \mathbb{P}\{\pi_t(\mathcal{Y}_{t-1};x) \neq \pi_t^*(x)|X_t = x\} P_X(dx)\Bigg]$$

$$\ge \sum_{t=1}^n \eta_t[d_t(\pi,\pi^*) - P_X\{[\theta-\eta_t,\theta+\eta_t]\}]$$

$$\ge \sum_{t=1}^n \eta_t[d_t(\pi,\pi^*) - C_*\eta_t^\alpha],$$

where the last inequality follows from (2.3).

Now we set $\varkappa = \max\{2, 1/(C_* x_0^\alpha)\}$, and $\eta_t = [d_t(\pi,\pi^*)/(\varkappa C_*)]^{1/\alpha}$. With this choice

$$\eta_t \le x_0 d_t^{1/\alpha}(\pi,\pi^*) \le x_0 \qquad \forall t \quad \text{and} \quad C_*\eta_t^\alpha \le \tfrac{1}{2} d_t(\pi,\pi^*).$$

Then it follows from (5.20) that

$$R_n(\pi,\pi^*) \ge \frac{1}{2}\sum_{t=1}^n \eta_t d_t(\pi,\pi^*) = \frac{n}{2}\left(\frac{1}{\varkappa C_*}\right)^{1/\alpha} \frac{1}{n}\sum_{t=1}^n [d_t(\pi,\pi^*)]^{1+1/\alpha}$$

$$\ge \frac{1}{2}\left(\frac{1}{\varkappa C_*}\right)^{1/\alpha} n^{-1/\alpha}[S_n(\pi,\pi^*)]^{1+1/\alpha},$$

where the last line follows from Jensen's inequality. □

*Proof of Theorem 3.* Fix $\delta > 0$, and let $\theta^{(0)} = 0$ and $\theta^{(1)} = \delta$. Note that when $\theta = \theta^{(0)}$ ($\theta = \theta^{(1)}$) it is preferable to sample from the arm 1 when $x > 0$ ($x > \delta$). Thus, $\pi^*(\theta^{(0)}, x) \neq \pi^*(\theta^{(1)}, x)$ only when $x \in (0,\delta)$.

Choose the probability density $f_{X,\delta}$ of $X$ so that

$$f_{X,\delta}(x) = \begin{cases} \frac{1}{2}C_*\alpha|x|^{\alpha-1}, & x \in [-x_0, \delta/2], \\ \frac{1}{2}C_*\alpha|x-\delta|^{\alpha-1}, & x \in [\delta/2, \delta + x_0]. \end{cases}$$

Suppose also that $\delta$ is small enough so that $C_* x_0^\alpha + C_*(\delta/2)^\alpha < 1$; then $f_{X,\delta}$ can be indeed continued outside the interval $[-x_0, x_0+\delta]$ so that it is a probability density. Clearly, the joint distributions $P_{X,Y}$ of $X$ and $Y$, corresponding to $\theta = \theta_{(0)}$ and $\theta = \theta_{(1)}$, and $f_{X,\delta}$ belong to $\mathcal{P}_\alpha$.

Therefore, we have

$$S_n(\pi;\mathcal{P}_\alpha) \ge \sup_{\theta\in\{\theta^{(0)},\theta^{(1)}\}} \sum_{t=1}^n \mathbb{P}_\theta\{\pi_t \neq \pi_t^*\}$$

$$\ge \frac{1}{2}\sum_{t=1}^n [\mathbb{P}_{\theta^{(0)}}\{\pi_t \neq \pi_t^*\} + \mathbb{P}_{\theta^{(1)}}\{\pi_t \neq \pi_t^*\}]$$

$$\ge \frac{1}{2}\sum_{t=1}^n \int_0^\delta [\mathbb{P}_{\theta^{(0)}}\{\pi_t(\mathcal{Y}_{t-1};x) \neq \pi_t^*(x)|X_t = x\}$$

$$+ \mathbb{P}_{\theta^{(1)}}\{\pi_t(\mathcal{Y}_{t-1};x) \neq \pi_t^*(x)|X_t = x\}] f_{X,\delta}(x)\,dx,$$

where we have used the fact that $X_t$ is independent of $\mathcal{Y}_{t-1}$. Here and in the sequel, $\mathbb{P}_{\theta^{(i)}}$ denotes the probability measure w.r.t. the distribution of observations $\mathcal{Y}_{t-1}$ when $\theta = \theta^{(i)}$.

Now for fixed $t$ consider the problem of testing the hypothesis $H_0 : \theta = \theta^{(0)}$ versus $H_1 : \theta = \theta^{(1)}$ from the observations $\mathcal{Y}_{t-1}$ collected under the policy $\pi$. Consider the following test: given observations $\mathcal{Y}_{t-1}$, the statistic $\pi_t(\mathcal{Y}_{t-1}; x)$ is computed for given $x \in (0, \delta)$, and it is compared with $\pi_t^*(\theta^{(0)}, x)$ and $\pi_t^*(\theta^{(1)}, x)$. The hypothesis $H_0$ is rejected when $\pi_t(\mathcal{Y}_{t-1}; x) \neq \pi_t^*(\theta^{(0)}, x)$. Because $\pi_t^*(\theta^{(0)}, x) \neq \pi_t^*(\theta^{(1)}, x)$, $\forall x \in (0, \delta)$, the expression under the integral sign in the last displayed formula above represents the sum of the error probabilities of the described test. Using well-known inequalities on error probabilities in testing problems [see, e.g., Devroye (1987) or Tsybakov (2004b)], we obtain that for any fixed $x \in (0, \delta)$

$$\begin{aligned}
&\mathbb{P}_{\theta^{(0)}}\{\pi_t(\mathcal{Y}_{t-1}; x) \neq \pi_t^*(x) | X_t = x\} + \mathbb{P}_{\theta^{(1)}}\{\pi_t(\mathcal{Y}_{t-1}; x) \neq \pi_t^*(x) | X_t = x\} \\
&\qquad \geq \tfrac{1}{4} \exp\{-\mathcal{K}\{\mathbb{P}_{\theta^{(0)}}(\mathcal{Y}_{t-1}), \mathbb{P}_{\theta^{(1)}}(\mathcal{Y}_{t-1})\}\},
\end{aligned}$$

where $\mathcal{K}\{\cdot, \cdot\}$ is the Kullback–Leibler divergence between distributions of the observations $\mathcal{Y}_{t-1}$ under $H_0$ and $H_1$. A straightforward calculation shows that

$$\begin{aligned}
&\mathcal{K}(\mathbb{P}_{\theta^{(0)}}(\mathcal{Y}_{t-1}), \mathbb{P}_{\theta^{(1)}}(\mathcal{Y}_{t-1})) \\
&\qquad = \mathbb{E}_{\theta^{(0)}}\left\{-\frac{1}{2\sigma^2}\sum_{s=1}^{t-1} \pi_s (Y_s - X_s)^2 + \frac{1}{2\sigma^2}\sum_{s=1}^{t-1} \pi_s (Y_s - X_s - \delta)^2\right\} \\
&\qquad \leq \frac{\delta^2}{2\sigma^2} \mathbb{E}_{\theta^{(0)}}[T_\pi(t-1)],
\end{aligned}$$

so that

$$\begin{aligned}
S_n(\pi; \mathcal{P}_\alpha) &\geq \frac{1}{8}\sum_{t=1}^{n} \exp\{-\mathcal{K}(\mathbb{P}_{\theta^{(0)}}(\mathcal{Y}_{t-1}), \mathbb{P}_{\theta^{(1)}}(\mathcal{Y}_{t-1}))\} \int_0^\delta f_{X,\delta}(x)\, dx \\
&\geq \frac{1}{8} C_*(\delta/2)^\alpha \sum_{t=1}^{n} \exp\left\{-\frac{\delta^2}{2\sigma^2}\mathbb{E}_{\theta^{(0)}} T_\pi(t-1)\right\} \\
&\geq \frac{1}{8} C_*(\delta/2)^\alpha n \exp\left\{-\frac{\delta^2 n}{2\sigma^2}\right\}.
\end{aligned}$$

Maximizing the RHS with respect to $\delta$, we set $\delta = \delta_* = \sqrt{\alpha}\sigma n^{-1/2}$. This yields

$$S_n(\pi; \mathcal{P}_\alpha) \geq \frac{1}{8}\left(\frac{\alpha}{2e}\right)^{\alpha/2} C_* \sigma^\alpha n^{1-\alpha/2},$$

as was claimed. The lower bound on $R_n(\pi; \mathcal{P}'_\alpha)$ follows from Lemma 5.

5.5. *Proof of Theorem 4.* $1^0$. We start with the proof of the lower bound on the regret.

Let $\pi$ be an arbitrary policy from $\Pi$. Without loss of generality, we assume that $x_0 = 1/2$ in the definition of the class $\mathcal{P}_1(\theta)$. For every fixed $\theta \in \Theta$ let $X_t$ be a random variable uniformly distributed on $A := [\theta - 1/2, \theta + 1/2]$, that is, $f_{X,\theta}(x) = I_A(x)$. Clearly, the corresponding joint distribution $P_{X,Y}$ of $(X, Y)$ belongs to the class $\mathcal{P}_1(\theta)$ (see Definition 1). For any fixed $\theta$ and $f_X = f_{X,\theta}$, we have

$$\begin{aligned} R_n(\pi, \pi^*) &= \mathbb{E}\sum_{t=1}^n |X_t - \theta| I\{\pi_t \neq \pi_t^*\} \\ &= \mathbb{E}\sum_{t=1}^n |X_t - \theta| I\{X_t \in [\hat{\gamma}_t \wedge \theta, \hat{\gamma}_t \vee \theta]\} \\ &= \mathbb{E}\sum_{t=1}^n \int_{\hat{\gamma}_t\wedge\theta}^{\hat{\gamma}_t\vee\theta} |x - \theta| I_A(x)\, dx = \frac{1}{2}\mathbb{E}\sum_{t=1}^n |\tilde{\gamma}_t - \theta|^2, \end{aligned}$$

where $\tilde{\gamma}_t = \max\{\hat{\gamma}_t, \theta + 1/2\}$ or $\tilde{\gamma}_t = \min\{\hat{\gamma}_t, \theta - 1/2\}$, depending on either $\hat{\gamma}_t \geq \theta$ or $\hat{\gamma}_t < \theta$. Thus, the problem is reduced to establishing a lower bound on the maximal cumulative squared error for estimating parameter $\theta \in \Theta$.

Let $\Lambda$ be a probability distribution on $\Theta$ with density $\lambda$ w.r.t. the Lebesgue measure. We assume that $\lambda$ converges to zero at the endpoints on the interval $\Theta$, and the Fisher information $I(\lambda)$ for the location parameter in $\lambda$ is positive and finite. The minimax cumulative risk in estimating $\theta$ is lower bounded by the Bayesian risk as follows:

$$\inf_{\tilde{\gamma}} \sup_{\theta\in\Theta} \mathbb{E}\sum_{t=1}^n |\tilde{\gamma}_t - \theta|^2 \geq \inf_{\tilde{\gamma}} \int \mathbb{E}\sum_{t=1}^n |\tilde{\gamma}_t - \theta|^2 \lambda(\theta)\, d\theta, \tag{5.21}$$

where inf is taken over all sequences $\tilde{\gamma} = (\tilde{\gamma}_t, t \geq 1)$ such that $\tilde{\gamma}_t$ is $\mathcal{F}_{t-1}$-measurable. Let $\mathcal{F}^*_{t-1} = \sigma(X_1, \ldots, X_{t-1}, Y_1, \ldots, Y_{t-1})$; because $\mathcal{F}_{t-1} \subset \mathcal{F}^*_{t-1}$, the expression on the RHS of (5.21) is lower bounded by $\inf_\gamma \int \mathbb{E}\sum_{t=1}^n |\gamma_t - \theta|^2$, where inf is taken over all sequences $\gamma = (\gamma_t)$ such that $\gamma_t$ is $\mathcal{F}^*_{t-1}$-measurable. Thus, we have

$$R_n(\pi; \mathcal{P}_1(\Theta)) \geq \sum_{t=1}^n \inf_{\gamma_t} \int \mathbb{E}|\gamma_t - \theta|^2 \lambda(\theta)\, d\theta, \tag{5.22}$$

where $\gamma_t$ is $\mathcal{F}^*_{t-1}$-measurable. Thus, the problem is reduced to establishing a lower bound on the Bayesian risk in the problem of estimating the scalar parameter $\theta \in \Theta$ from observations $\{(X_s, Y_s), s = 1, \ldots, t-1\}$, where $Y_s = X_s - \theta + \varepsilon_s$, and $\varepsilon_s$ are i.i.d. zero mean Gaussian random variables with variance $\sigma^2$. This problem is well studied, and there are different methods for establishing such lower bounds [see, e.g., Borovkov and Sakhanenko (1980), Brown and Gajek (1990) and Gill and Levit (1995)].

In particular, by the van Trees inequality [see Gill and Levit (1995)],

$$\inf_{\gamma_t} \int \mathbb{E}|\gamma_t - \theta|^2 \lambda(\theta)\, d\theta \geq \frac{1}{I_t + I(\lambda)},$$

where $I_t$ is the expected Fisher information for $\theta$ associated with the conditional density of observations $(X_1, Y_1, \ldots, X_{t-1}, Y_{t-1})$ given $\theta$; and $I(\lambda)$ is the Fisher information for the location parameter in $\lambda$. Thus,

$$I_t := \mathbb{E}\left[\frac{\partial}{\partial\theta} \ln f(\mathcal{Y}_{t-1}|\theta)\right]^2 = \mathbb{E}\left[-\frac{1}{\sigma^2}\sum_{s=1}^{t-1}(Y_s - X_s - \theta)\right]^2 = \frac{t-1}{\sigma^2}.$$

The standard choice of $\lambda$ is the following:

$$\lambda(\theta) = \frac{1}{h}\lambda_0\left(\frac{\theta - \tau_0}{h}\right), \qquad \lambda_0(\theta) = \cos^2(\pi\theta/2) I\{|\theta| \leq 1\},$$

where $\tau_0$ is the center of the interval $\Theta := [\tau^-, \tau^+]$, and $h = \tau^+ - \tau^-$. With this choice $I(\lambda_0) = \pi^2$ and $I(\lambda) = h^{-2} I(\lambda_0) = \pi^2 h^{-2}$. Therefore, applying the van Trees inequality for each summand in (5.22), we obtain

$$R_n(\pi; \mathcal{P}_1(\Theta)) \geq \sum_{t=1}^{n} \frac{1}{I_t + I(\lambda)} = \sigma^2 \sum_{t=2}^{n} \frac{1}{t - 1 + \sigma^2\pi^2 h^{-2}} \geq c\sigma^2 \ln n$$

for $n$ large enough.

$2^0$. The lower bound on $S_n(\pi; \mathcal{P}_2(\Theta))$ follows from identical considerations. In this case we choose $f_X$ to be linear in the vicinity of $\theta$; then for any policy $\pi \in \Pi$, $S_n(\pi, \pi^*) \geq c\mathbb{E}\sum_{t=1}^{n} |\tilde{\gamma}_t - \theta|^2$ for any sequence $\tilde{\gamma} = (\tilde{\gamma}_t)$ of random variables such that $\tilde{\gamma}_t$ is $\mathcal{F}_{t-1}$-measurable.

Department of Statistics
University of Haifa
Haifa 31905
Israel
E-mail: goldensh@stat.haifa.ac.il

Graduate School of Business
Columbia University
New York, New York 10027
USA
E-mail: assaf@gsb.columbia.edu